\def\ZE{{\mathbb Z}}
\def\RE{{\mathbb R}}
\def\CX{{\mathbb C}}
\def\pol{{\textstyle\frac{1}{2}}}

\def\CP{\mathbb{C}P}
\def\HC{\mathbb{H}}
\def\ddb{\frac{i}{2\pi}\partial\bar{\partial}}
\def\tddb{\textstyle{\frac{i}{2\pi}\partial\bar{\partial}}}
\def\pbp{\partial\bar{\partial}}

\def\hra{\hookrightarrow}
\def\eps{\varepsilon}
\def\BE{\eus{E}}
\def\BF{\eus{F}}

\def\fcy{f_{\mathrm{cyl}}}
\def\omc{\omega_{\mathrm{cyl}}}
\def\omca{\omega_{\mathrm{cyl}}(a)}
\def\omco{\omega_{\mathrm{cyl},0}}
\def\cpt{\omega_{\mathrm{comp}}}
\def\oW{\overline{W}}
\def\tW{W'}
\def\p{\partial}
\def\minus{\setminus}

\def\tD{{\tilde{D}}}
\def\tV{{\tilde{V}}}
\def\T{\tilde{\Theta}}
\def\hra{\hookrightarrow}
\def\xra{\xrightarrow}

\documentclass[12pt]{amsart}
\usepackage[english]{babel}
\usepackage{a4wide,amscd,amssymb,array,euscript}

\DeclareFontFamily{OT1}{rsfs}{}
\DeclareFontShape{OT1}{rsfs}{n}{it}{<-> rsfs10}{}
\DeclareMathAlphabet{\curly}{OT1}{rsfs}{n}{it}

\def\eus{\EuScript}
\newtheorem{thrm}[equation]{Theorem}
\newtheorem{lemma}[equation]{Lemma}
\newtheorem{prop}[equation]{Proposition}
\newtheorem{cor}[equation]{Corollary}
\theoremstyle{definition}
\newtheorem*{defi}{Definition}
\newtheorem{example}[equation]{Example}
\theoremstyle{remark}
\newtheorem*{remark}{Remark}
\newtheorem*{remarks}{Remarks}
\let\phi=\varphi
\let\epsilon=\varepsilon

\DeclareMathOperator{\const}{const}
\DeclareMathOperator{\Hol}{Hol}
\DeclareMathOperator{\Ker}{Ker}
\DeclareMathOperator{\im}{Im}
\DeclareMathOperator{\re}{Re}

\DeclareMathOperator{\rank}{rk}
\DeclareMathOperator{\Ric}{Ric}
\DeclareMathOperator{\SP}{Sp}

\DeclareMathOperator{\vol}{vol}
\newcommand{\HH}{\eus{H}}
\newcommand{\sieq}{\lower-.25ex\hbox{$\simeq$}}
\def\K3{\curly{K}\!\text{\it 3}}
\newcommand{\MM}{\curly{P}(\curly{V})}
\newcommand{\Mm}{\curly{P}(\curly{V}_1)}
\def\OH{\curly{O}}
\def\U{\curly{U}}
\def\V{\curly{V}}
\begin{document}

\title{Twisted connected sums and special Riemannian holonomy}
\author[A.~Kovalev]{Alexei Kovalev}
\address{Department of Pure Mathematics and Mathematical Statistics,
	Centre for \mbox{Mathematical} Sciences, Wilberforce Road,
	Cambridge CB3 0WB, UK}
\email{A.G.Kovalev@dpmms.cam.ac.uk}

\begin{abstract}
We give a new, connected sum construction of Riemannian metrics with
special holonomy $G_2$ on compact 7-manifolds.  The construction is
based on a gluing theorem for appropriate elliptic partial
differential equations.  As a prerequisite, we also obtain
asymptotically cylindrical Riemannian manifolds with holonomy $SU(3)$
building up on the work of Tian and Yau.  Examples of new topological
types of compact 7-manifolds with holonomy $G_2$ are constructed using
Fano 3-folds.
\end{abstract}

\maketitle

The purpose of this paper is to give a new construction of compact
7-dimensional Riemannian manifolds with holonomy group~$G_2$.  The
holonomy group of a Riemannian manifold is the group of isometries of
a tangent space generated by parallel transport using the Levi--Civita
connection over closed paths based at a point. For an oriented
$n$-dimensional manifold the holonomy group may be identified as a
subgroup of $SO(n)$. If there is a structure on a manifold defined by
a tensor field and parallel with respect to the Levi--Civita
connection then the holonomy may be a proper subgroup of $SO(n)$ (it
is just the subgroup leaving invariant the corresponding tensor
on~$\RE^n$). There is essentially just one possibility for such
holonomy reduction in odd dimensions, as follows from the well-known
Berger classification theorem. This `special holonomy' group is the
exceptional Lie group $G_2$ and it occurs in dimension $n=7$, when a
metric is `compatible' with the non-degenerate cross-product
on~$\RE^7$, see \S\ref{synopsis} for the precise definitions.  The
only previously known examples of compact Riemannian
\mbox{7-manifolds} with holonomy~$G_2$ are due to Joyce who used a
generalized Kummer construction and resolution of
singularities~\cite{joyce,joyce-book} (for non-compact examples
see e.g.~\cite{baer,GG,bryant,bryant-salamon,gibbons}).
The compact 7-manifolds with holonomy $G_2$ in this paper are obtained
by a different, connected-sum-like construction for a pair of
non-compact manifolds with asymptotically cylindrical ends.

It is by now understood that the existence problem for the metrics
with holonomy~$G_2$ can be expressed as a non-linear system of PDEs
for a non-degenerate differential 3-form $\phi$ on a 7-manifold.
More precisely, a solution 3-form defines a torsion-free
$G_2$-structure and thus a metric, from the inclusion $G_2\subset
SO(7)$. The holonomy of this metric is in general only {\em contained}
in~$G_2$.  In fact, we shall obtain the solutions, torsion-free
$G_2$-structures, by proving a gluing theorem for pairs of manifolds
with holonomy $SU(3)$, a maximal subgroup of~$G_2$. To claim that the
holonomy on the resulting compact 7-manifold is exactly~$G_2$ it
suffices to verify a topological condition: that the 7-manifold has
finite fundamental group.\pagebreak

To implement a connected sum strategy we require, in the first place,
a suitable class of non-compact, asymptotically cylindrical Riemannian
manifolds with holonomy $SU(3)$. The point is that while the task of
constructing asymptotically cylindrical metrics with holonomy $G_2$ in
general is not likely to be any easier than the search for holonomy
$G_2$ on compact manifolds, the metrics with holonomy $SU(3)$ are
understood much better. The special holonomy $SU(3)$ naturally occurs
on manifolds of {\em complex} dimension~3, for the Ricci-flat K\"ahler
metrics. In~\cite{TY1,TY2} Tian and Yau obtained a number of
existence theorems for complete Ricci-flat K\"ahler metrics on
quasiprojective manifolds.  With some additional work, we are able to
find simple topological conditions on these manifolds to further
ensure that we obtain Ricci-flat K\"ahler metrics with holonomy
$SU(3)$ and with the desired cylindrical asymptotic model at infinity.
The result, Theorem~\ref{ric-flat}, produces a non-compact analogue of
Calabi--Yau 3-folds and may be of independent interest---from the
point of view of Riemannian geometry and analysis, these are examples
of manifolds with `added boundary at infinity', or $b$-manifolds in
the sense of Melrose~\cite{tapsit}.

The manifolds with holonomy $SU(3)$ that we obtain have real
dimension~6. The Riemannian product with a circle then yields
7-manifolds with the same holonomy and with an end asymptotic to
a half-cylinder. The cross-section of this half-cylinder is the
Riemannian product of two circles and a complex surface of type K3
with a hyper-K\"ahler metric. A closed compact Riemannian 7-manifold
may be constructed from a pair of these asymptotically cylindrical
7-manifolds by truncating the ends, cutting off to the cylindrical
metric and identifying the boundaries via an orientation-reversing
isometry. Thus the connected sum operation that we perform is not a
usual one as the cross-section of the neck is not a sphere. The second
and more important distinction is the non-trivial, `twisting' isometry
map identifying the two boundaries. The basic idea is that the map
interchanges the two $S^1$ factors to avoid an infinite fundamental
group of the compact 7-manifold. The map also interchanges the complex
structures (using Torelli theorem) on the two K3 surfaces in such a
way that the compact 7-manifold has a well-defined $G_2$-structure
compatible with that on the summands. The precise construction can be
found in~\S\ref{K3}.

The $G_2$-structure thus constructed on the compact 7-manifold
can be regarded as an `approximate solution' to the PDEs for
torsion-free $G_2$-structures. The setting achieved by the above steps
falls under a general pattern of gluing the solutions of elliptic PDEs
for special Riemannian structures on (generalized) connected sums.
The gluing argument is carried out in~\S\ref{deform}, using the
analytical techniques developed in~\cite{asdsums} and previously in
a special case in~\cite{floer-conformal}. 
It is worth to point out that, with a careful choice of the Banach
spaces of differential forms, our gluing is unobstructed. Therefore,
we can be sure to solve the equation and obtain a metric with holonomy
$G_2$ on our compact 7-manifold if only the initial approximate
solution is well-defined.

The examples we found for the connected sum construction are discussed
in general in~\S\ref{fano} and some specific calculations are
performed in~\S\ref{examples}. They originate from the class of compact
complex 3-folds with Ricci-positive K\"ahler metric, known as Fano
manifolds. Fano manifolds were extensively studied over the past
decades in relation to problems in algebraic geometry and
K\"ahler--Einstein metrics. Exploiting the
theory of Fano manifolds and K3 surfaces, and also an appropriate
extension of the deformation theory of complex manifolds, we are able
to prove a general Theorem~\ref{ult} providing the required initial
setting for the `twisted connected sum' from any pair of algebraic
families of smooth Fano 3-folds.
This quite rapidly leads to many new topological types
of compact manifolds of holonomy $G_2$.  The geometrical consequences
concerning a distinguished class of 
minimal submanifolds, called coassociative
submanifolds, of these new 7-manifolds will be investigated in a
sequel~\cite{inprep} to this paper.

\medskip

{\bf Acknowledgements.} I am indebted to Simon Donaldson for
explaining to me that the interchanging of the circles and
complex structures in the connected sum would lead
to a construction of $G_2$-manifolds and leaving it to me to further
develop the idea. I am also grateful to Dominic Joyce for discussions
on his work, to Elmer Rees for conversations leading to several
improvements, to Alessio Corti for a useful suggestion concerning
the proof of Theorem~\ref{onto}, and to Mark Gross for enlightening
remarks. I thank the referee for  suggesting several corrections.
The author was supported by the EPSRC for part of
the duration of this work. Partially supported by the European
Contract Human Potential Programme, Research Training Network
HPRN-CT-2000-00101.

\section{Torsion-free $G_2$-structures}
\label{synopsis}

This section gives a short summary of some standard results on the
Riemannian geometry associated with the group $G_2$. Further details
can be found in~\cite{calibr}, \cite{joyce-book} and~\cite{salamon}.

The group $G_2$ may be defined as the group of linear automorphisms
of the cross-product algebra on~$\RE^7$, the vectors in~$\RE^7$ being
interpreted as pure imaginary octonions. Any automorphism in $G_2$
necessarily preserves the Euclidean metric and so $G_2$ is
precisely the group leaving invariant the differential 3-form $\phi_0$
on $\RE^7$,
$$
\langle a\times b, c\rangle = \phi_0(a,b,c),
$$
encoding the cross-product multiplication.
The Euclidean metric is completely determined by $\phi_0$, in
particular, there is an explicit relation
	\begin{equation}
	\label{metr}
6\langle a, b \rangle\, d\!\vol
 = (a\lrcorner\phi_0)\wedge(b\lrcorner\phi_0)\wedge\phi_0.
	\end{equation}
The group~$G_2$ is a closed proper (14-dimensional)
subgroup of $SO(7)$ and thus a compact Lie group.

The Hodge star on $\RE^7$ is $G_2$-invariant, thus $G_2$ is also
obtained as the stabilizer
of the 4-form $\sigma_0=*\phi_0$ on~$\RE^7$ in the $GL_+(7,\RE)$-action.
The $GL_+(7,\RE)$-orbits of both $\phi_0$ and $\sigma_0$ are open in
the space of all 3- and 4-forms on $\RE^7$ respectively.
($GL_+(7,\RE)\subset GL(7,\RE)$ denotes the subgroup of
orientation-preserving linear isomorphisms.)

If $M$ is an oriented 7-manifold then a smoothly varying non-degenerate
cross-product on the tangent spaces defines a $G_2$-structure on~$M$.
Equivalently, a $G_2$-structure on~$M$ may be given by a 3-form $\phi$
whose value $\phi_p$ at each point $p\in M$ may be written as $\phi_0$
with respect to some positively oriented tangent frame giving an
isomorphism $T_pM\cong\RE^7$. The space of
all $G_2$-structures on $M$ therefore may be identified with a
subset of differential 3-forms, equivalent to $\phi_0$ at every
point of~$M$, and we shall denote this subset by
$\Omega_+^3(M)\subset\Omega^3 (M)$ (it is not a linear subspace).
\pagebreak

Any $G_2$-structure, being a special case of $SO(7)$-structure,
induces a metric on~$M$ (the formula~\eqref{metr} recovers
this metric up to a conformal factor).
Given a $G_2$-structure $\phi$ on~$M$, $\phi+\chi$ is again a
$G_2$-structure whenever $\|\chi\|_{g(\phi)}<\epsilon$ in the sup-norm in
the metric induced by~$\phi$. In this sense, $\Omega_+^3(M)$ is an open
subset of $\Omega^3 (M)$. The constant $\epsilon$ is determined by 
a calculation in $\RE^7$ and is independent of $\phi$ or~$M$.

There is a completely equivalent, dual description of
$G_2$-structures on~$M$ as 4-forms in $\Omega_+^4(M)\subset\Omega^4 (M)$
modelled point-wise on $\sigma_0$. That is, a 4-form is in
$\Omega_+^4(M)$ precisely when it may be realized as
$\sigma_0$ at any point of~$M$, by choosing a
positively oriented tangent frame and thus an identification of
tangent space with~$\RE^7$.
Similarly to $\Omega_+^3(M)$ above, the set $\Omega_+^4(M)$ is
open in $\Omega^4 (M)$ in the sup-norm.

The equivalence of the descriptions of $G_2$-structures as 3-forms and
4-forms is given by the map
$$
\Theta: \phi\in\Omega_+^3 \to *_{g(\phi)}\phi\in\Omega_+^4,
$$
which is non-linear because of the dependence of the Hodge star on
$\phi$ through the induced metric $g(\phi)$.

A metric $g(\phi)$ coming from a $G_2$-structure $\phi$ does not
necessarily have holonomy contained in $G_2$. For that to be the case,
the $G_2$-structure $\phi$ must be
{\em torsion-free}~\cite{salamon}, Lemma~11.5, a condition expressed
by the differential equations
	\begin{equation}
	\label{torsion}
d\phi = 0 \qquad\text{and}\qquad d\Theta(\phi) = 0.
	\end{equation}
For a compact 7-manifold $M$ with a torsion-free
$G_2$-structure~$\phi$, the holonomy $\operatorname{Hol}(g(\phi))$ is
exactly $G_2$ if and only if $\pi_1(M)$ is finite~\cite{joyce},
II~Proposition~1.1.1.
We use the terms \mbox{\em $G_2$-manifold} and \mbox{\em $G_2$-metric}
in this paper in the `strong' sense to mean a Riemannian manifold
$(M,g)$ such that $\operatorname{Hol}(g)=G_2$ (and likewise for other
holonomy groups).

\section{Asymptotically cylindrical manifolds with holonomy $SU(3)$}
\label{kahler}

If a torsion-free $G_2$-structure $\phi$ is given on a product of a
6-manifold and a circle $S^1$ and induces a product metric then
there is a parallel non-zero vector field tangent to the $S^1$
direction. Respectively, the holonomy of the metric on the 6-manifold
is contained in a maximal subgroup $SU(3)\subset G_2$ of the
automorphisms of the cross-product algebra $\RE^7$ fixing a non-zero
vector. In particular, the metrics with holonomy in $SU(3)$ on a
6-manifold give rise to solutions of~\eqref{torsion}.
In this section we construct a class of such metrics. 

A metric $g$ on a real $2n$-dimensional manifold $W$ has holonomy
in~$U(n)$ if and only if that metric is K\"ahler. In particular, $W$
then admits an integrable complex structure $I$ orthogonal with
respect to $g$ and a closed $(1,1)$-form, the K\"ahler form $\omega$,
parallel with respect to $g$.  Ricci curvature of a K\"ahler metric is
equivalent to the curvature of the induced Hermitian connection on the
canonical bundle $K_W$ of $(n,0)$-forms, so $g$ is Ricci-flat
precisely when $K_W$ is a flat Hermitian line bundle.  The sections of
$K_W$ parallel with respect to $g$ are holomorphic, hence closed,
forms.  If there is a globally defined parallel
$(n,0)$-form $\Omega$ on~$W$, sometimes called a
`holomorphic volume form', then the holonomy of~$g$ is contained in
$SU(n)$.  This will be the case e.g.\ when $W$ is simply-connected.
Then $\Omega\wedge\Omega^*$ is a $g$-parallel $(n,n)$-form, hence a
constant multiple of the volume form $\omega^n$, where $\Omega^*$
denotes the $(0,n)$-form complex conjugate of~$\Omega$.

The $SU(3)$-metrics obtained in this section are asymptotically
cylindrical Ricci-flat K\"ahler metrics on complex 3-dimensional
manifolds. 
These manifolds are non-compact and may be given as $W=\oW\minus D$,
where $\oW$ is a compact K\"ahler manifold containing a compact
complex surface~$D$. We shall make an additional assumption that $D$
can be given as the zero set of a holomorphic coordinate $z$ on
$\oW$ and a tubular neighbourhood $U$ of $D$ in $\oW$ is diffeomorphic,
as a real manifold, to $D\times \{|z|<1\}$. With a change of
variable $z=e^{-t-i\theta}$, $t>0$, $\theta\in S^1$, $W$
can be considered as a manifold with cylindrical end $U\setminus D$
asymptotic to $D\times\RE_{>0}\times S^1$.  We
want to put on $W$ a Ricci-flat K\"ahler metric asymptotic
to a product Ricci-flat K\"ahler metric on the half-cylinder
$D\times\RE_{>0}\times S^1$.

	\begin{subequations}\label{su2su3}
The holonomy representation for a product metric on
$D\times\RE_{>0}\times S^1$ reduces to the holonomy representation
for a metric on~$D$. We shall be interested in the case when $D$ 
is a simply-connected compact complex surface with $c_1(D)=0$, that
is, a K3 surface. It is a well-known consequence of Yau's solution of
the Calabi conjecture~\cite{yau} that a K3 surface admits a unique
Ricci-flat K\"ahler metric in every K\"ahler class. We shall write
$\kappa_I$ for the K\"ahler form
and $\kappa_J+i\kappa_K$ for a holomorphic volume $(2,0)$-form of this
metric, where $\kappa_I,\kappa_J,\kappa_K$ are real 2-forms on~$D$.
The subscript notation indicates the fact that a Ricci-flat K\"ahler
metric on a K3 surface is hyper-K\"ahler, compatible with the action of
quaternions on tangent spaces, in particular one has
$\kappa_I^2 = \kappa_J^2 = \kappa_K^2$. We shall consider the
hyper-K\"ahler geometry of K3 surfaces in more detail in~\S\ref{K3}.
For the moment, we only note that a product Ricci-flat K\"ahler metric
on the half-cylinder $\RE_{>0}\times S^1\times D$ has K\"ahler form
	\begin{equation}
\omega_0=dt\wedge d\theta + \kappa_I,
	\end{equation}
and holomorphic volume form
	\begin{equation}
\Omega_0=(dt+id\theta)\wedge(\kappa_J+i\kappa_K).
	\end{equation}
	\end{subequations}

The following theorem provides a class of manifolds carrying the
$SU(3)$-structures modelled on $(\omega_0,\Omega_0)$ `near infinity'.
It builds up on the work of Tian and Yau~\cite{TY1,TY2} on the
existence of complete K\"ahler metrics with prescribed Ricci curvature.
We shall sometimes refer to K\"ahler metrics by their K\"ahler forms.
The notation $O(e^{-\lambda t})$ means that a function or differential
form is bounded by $\const\cdot\, e^{-\lambda t}$.
	\begin{thrm}\label{ric-flat}
	\begin{subequations}\label{atinfty}
Let $\oW$ be a smooth compact K\"ahler 3-fold with $H^1(\oW,\RE)=0$
and $\omega'$ the K\"ahler form on~$\oW$. Suppose that a K3 surface
$D$ in~$\oW$ is an anticanonical divisor and has trivial
self-intersection class $D\cdot D=0$ in $H_2(\oW,\ZE)$.

Then $W=\oW\minus D$ admits an asymptotically cylindrical Ricci-flat
K\"ahler metric~$g$ with respect to a diffeomorphism
$U\minus D\to D\times\{0<|z|<1\}$, where $z$ is a holomorphic function
on a neighbourhood $U$ of~$D$ in~$\oW$ and $z$ vanishes to order~1
on~$D$. The K\"ahler form of~$g$ can be written near $D$ as
	\begin{equation}\label{asymp}
\omega_g|_{U\minus D} = \omega_0 + d\psi.
	\end{equation}
If $g$ has a holomorphic volume form $\Omega_g$ then (up to a constant
factor)
	\begin{equation}\label{h.vol}
\Omega_g|_{U\minus D} = \Omega_0 + d\Psi.
	\end{equation}
Here $\omega_0,\Omega_0$ are identified with the
differential forms on the (half-)cylinder given by~\eqref{su2su3}
with $e^{-t-i\theta}=z$, and the K\"ahler class
$[\kappa_I]\in H^{1,1}(D)$ is~$[\omega'|_D]$.
The 1-form $\psi$ and 2-form $\Psi$ are
smooth and $O(e^{-\lambda t})$ with all derivatives (as measured by
the metric~$\omega_0$), for any $\lambda< \min\{1,\sqrt{\lambda_1(D)}\}$,
where $\lambda_1(D)$ is the first
positive eigenvalue of the Laplacian on~$D$ with the metric~$\kappa_I$.
	\end{subequations}
	\end{thrm}
Theorem~\ref{ric-flat} is proved in the next section.
	\begin{remark}
There is a convenient way to formalize the idea of the `boundary at
infinity' of the Riemannian manifold $(W,g)$. A change of real
coordinate $x=e^{-t}$, extended to map $t=\infty$ to~\mbox{$x=0$,}
transforms $W$ into a compact manifold $\tilde{W}$ with boundary
$S^1\times D$ given as the zero set of~$x$. The asymptotically
cylindrical metric $g$ given by Theorem~\ref{ric-flat} can be written
near the boundary as
	\begin{equation}\label{bmetr}
\smash[t]{g=\bigl(\frac{dx}x\bigr)^2+\tilde{g}.}
	\end{equation}
Here $\tilde{g}$ is a symmetric semi-positive bilinear form on
$\tilde{W}$ which is $C^\infty$-smooth up to the boundary and
$\tilde{g}|_{x=0}$ is a product metric on $S^1\times D$,
with the K\"ahler metric $\kappa_I$ on~$D$. Elliptic theory on
manifolds with boundary endowed with Riemannian metrics of
type~\eqref{bmetr} is developed in~\cite{tapsit}, where these are
called `exact $b$-metrics'. Results proved for manifolds with
exact $b$-metrics are therefore valid, after a change of notation,
on the asymptotically cylindrical manifolds obtained in
Theorem~\ref{ric-flat}. This will become important several times
for the analysis in \S\ref{bCY} and~\S\ref{deform} below.
	\end{remark}
Examples of manifolds $W$ satisfying the hypotheses of
Theorem~\ref{ric-flat} will be given in~\S\ref{fano}. To ensure
that they have a well-defined $SU(3)$-structure $(\omega_g,\Omega_g)$
we shall need the following.
	\begin{thrm}
        \label{line}
Suppose that a 3-fold $\oW$ and K3 surface $D$ in~$\oW$ satisfy
the hypotheses of Theorem~\ref{ric-flat} and let $g$ be
the asymptotically cylindrical Ricci-flat K\"ahler metric on
$W=\oW\minus D$. If, in addition, $\oW$ is simply-connected and
contains a complex curve $\ell$ with $D\cdot \ell=m>0$ then the metric
$g$ has holonomy~$SU(3)$.
	\end{thrm}
	\begin{pf}
Consider $\oW$ as a union $\oW=W\cup U$ where $U$ is a small tubular
neighbourhood of~$D$. By the Seifert and van Kampen
theorem (e.g.~\cite{armstrong}, Theorem 6.13), the fundamental group
of $\oW$ may be expressed as the `amalgamated product',
$\pi_1(\oW)=\pi_1(W)*_{\pi_1(W\cap U)}\pi_1(U)$, over $\pi_1(W\cap U)$.
Since $\pi_1(U)$ is trivial, the amalgamated product in this case is
just the quotient of $\pi_1(W)$ by the image of $\pi_1(W\cap U)$ under
homomorphism induced by the inclusion $W\cap U\subset W$. By the
hypothesis, this quotient is the trivial $\pi_1(\oW)$. Hence
$\pi_1(W)$ is isomorphic to the quotient of $\pi_1(W\cap U)=\ZE$ by a
subgroup, so $\pi_1(W)$ is a {\em cyclic} group, generated by a loop
around~$D$ in~$\oW$. Then $\pi_1(W)\cong H_1(W)$ with the coefficients
in~$\ZE$.  Furthermore, this cyclic group is {\em finite} as
may be seen by considering a part of Mayer--Vietoris exact sequence
for $\oW=W\cup U$,
$$
H_2(\oW)\to H_1(U\minus D) \to H_1(W)\oplus H_1(U)\to H_1(\oW).
$$
The first homomorphism maps to $H_1(U\minus D)=\ZE$ evaluating the
intersection number with $[D]$ and its image contains $m\ZE$ by
the hypothesis. To simplify the notation we shall pretend that the
image actually is $m\ZE$, as the image is a non-trivial subgroup
of~$\ZE$ anyway. As $H_1(U)=H_1(\oW)=0$, we must have $\pi_1(W)\cong
H_1(W)=\ZE_m$.

If $m=1$, so $W$ is simply-connected, then there is nothing to prove
as the vanishing Ricci curvature makes $\Hol(g)=SU(3)$ the only
possibility, by Berger classification of Riemannian holonomy
groups, \cite{salamon}, Ch.10. For $m>1$, we shall be done if we show
that the flat Hermitian connection induced by~$g$ on the canonical
bundle $K_W$ is a trivial, product connection. It suffices to consider
the holonomy transformation for this connection along a generator of
$\pi_1(W)$. The loop around~$D\subset\oW$ corresponds to the $S^1$
factor in the cross-section of the cylindrical end of~$W$. The
holonomy transformation along this $S^1$ factor equals $m$th root of
unity and is independent of the choice of loop in the homotopy class.
On the other hand, as $g$ decays to a product Ricci-flat metric on
$S^1\times D$, 
the same holonomy transformation tends to the identity when the loop
goes to the infinity along the end of~$W$. So, this transformation
must be the identity and hence the holonomy of $g$ is in $SU(3)$. Then
by considering the universal Riemannian covering of~$W$ we conclude
that $\Hol(g)$ is precisely~$SU(3)$.
	\end{pf}

\section{Proof of Theorem~\ref{ric-flat}}
\label{bCY}

Let $s\in H^0(\oW,K^{-1}_{\oW})$ denote the defining section for~$D$.
As the surface $D$ has trivial self-intersection class in~$H_2(\oW)$,
the normal bundle of~$D$ has degree zero, by Mumford's
self-intersection formula, \cite{H}, p.431.  Then, by the adjunction
formula, the anticanonical bundle 
$K^{-1}_{\oW}$ restricts to a holomorphically trivial line bundle
on~$D$.  With the help of a holomorphic
trivialization of $K^{-1}_{\oW}$ over a tubular neighbourhood~$U$
of~$D$, the section $s$ defines a holomorphic function, $\tau$ say,
on~$U$ with $D=\tau^{-1}(0)$. The vanishing of Dolbeault cohomology
$H^{0,1}(\oW)$ implies, cf.~\cite{GH} pp.34--35, that $1/\tau$ can be
extended to a meromorphic function on all of~$\oW$ with the only pole
along~$D$ of order~1. Thus $\tau$ extends to a K3 fibration of~$\oW$,
which we denote again by $\tau:\oW\to\CP^1$, with $D$ a fibre.
It also follows that $K^{-1}_{\oW}$ is actually the pull-back
via~$\tau$ of a holomorphic line bundle of degree~1 over $\CP^1$ and
$s=s_1\circ\tau$ for a section $s_1$ over~$\CP^1$.  We shall use $z$
interchangeably to denote a holomorphic coordinate on $\CX P^1$ and
the corresponding pull-back holomorphic coordinate on~$\oW$. We can
use the triviality of the normal bundle of~$D$ to define a local
product decomposition
	\begin{equation}\label{prod}
U\simeq\{|z|<1\}\times D,
	\end{equation}
identifying $\tau^{-1}(z)$ with $\{z\}\times D$; note that
\eqref{prod} is only a diffeomorphism of the underlying {\em real}
manifolds as the holomorphic structure of $\tau^{-1}(z)$ generally
depends on~$z$. (However $z$ is holomorphic in both complex
structures in~\eqref{prod}.) The metric $\omega_0$ defined
in~(\ref{su2su3}a) is K\"ahler with respect to the complex structure
of $\{|z|<1\}\times D$,  $z=e^{-t-i\theta}$.

Recall that if $\omega$ is a K\"ahler metric on a complex $n$-manifold
and $e^f\omega^n$ is the volume form of another K\"ahler metric then
the Ricci form of the latter metric is $\Ric(\omega)-i\p\bar{\p}f$. 
The following existence result is a direct application
of~\cite{TY1}, Theorems~1.1 and~5.2 (cf.\ also~\cite{TY2}, p.~52). 
	\begin{prop}\label{soln}
There is a choice of K\"ahler metric $\omc$ on $W$ such that
$\omc|_{U\minus D}$ is commensurate with the cylindrical
metric~$\omega_0$ and $\Ric(\omc)=i\p\bar{\p}f$ with
$\int_W(e^f-1)\omc^3=0$, where 
$|f|=O(t^{-N})$ for some $N>4$.

Further, the complex Monge--Amp\`ere equation
	\begin{equation}\label{MA}
(\omc+\ddb u)^3 = e^{f}\omc^3,
	\end{equation}
has a smooth solution $u$ on $W$, such that on
$U\minus D\;\sieq\;\RE_{>0}\times S^1\times D$ the function $u$
converges to zero uniformly in $S^1\times D$, as $t\to\infty$. The
derivatives of~$u$ are bounded (with respect to $\omc$) and
$\omc+\tddb u>0$. Thus $\omega_g=\omc+\ddb u$ is a Ricci-flat K\"ahler
metric on~$W$ commensurate with $\omc$.
        \end{prop}
	\begin{prop}\label{unique}
The solution $u$ is uniquely determined by $\omc$ and $f$, where
$\omc$, $f$, and $u$ are as defined in Proposition~\ref{soln}.
	\end{prop}
	\begin{pf}[Proof of Proposition~\ref{unique}]
Suppose that $u'$ is another solution of~\eqref{MA} and has the same
properties as $u$, as given in Proposition~\ref{soln}. Recall that
on a compact manifold the uniqueness of solution of~\eqref{MA}, up to
additive constant, is proved using local calculations and Stokes'
Theorem to argue that the
integral of $\Delta_g(u-u')^2$ is zero, where $\Delta_g$ denotes the
Laplacian of the metric~$\omega_g$. See \cite{yau}, Theorem~3, or
\cite{joyce-book}, \S 5.7. The argument of the compact case remains
valid in our setting because the solutions $u,u'$ were assumed to
decay to zero along the cylindrical end of~$W$.
	\end{pf}

We claim that $\omega_g$ is the metric that we want, i.e.\
\eqref{asymp} holds.  We shall first deduce, by examining the method
of~\cite{TY1}, that the rates of decay of $\omc$ and $f$ can be taken
to be exponential in~$t$. In the present situation we can
exploit the K3 fibration structure of~$\oW$ and the calculations are
simplified.

Consider on~$\oW$ a smooth function $u_0$ with support in~$U$ and such
that $(\omega' + \tddb u_0)|_D$ is the Ricci-flat K\"ahler metric
(Calabi--Yau metric) in the K\"ahler class $[\omega'|_D]$ on the K3
surface~$D$. The form $\omega'+\tddb u_0$ is positive on the fibres
of~$\tau$ contained in some tubular neighbourhood of~$D$. By smoothly
cutting off to zero within this tubular neighbourhood, but away
from~$D$, the function $u_0$ may be chosen so that the form
$\omega'+\tddb u_0$ is positive on any fibre of~$\tau$ in~$\oW$.
\enlargethispage{1.ex}

Denote by $\omega_1$ a K\"ahler form on~$\CP^1$, such that
$\omega_1=(1+O(|z|^2))idz\wedge d\bar{z}$ for small $|z|$, where
$\tau^{-1}(\{z=0\})=D$. The pull-back $\tau^*\omega_1$ is a semi-positive
(1,1)-form on~$\oW$, positive in the directions transverse to the
fibres of~$\tau$. Rescaling $\omega_1$ by an appropriate positive
constant $\mu$, we obtain a positive (1,1)-form 
$$
\cpt = \cpt(\mu) = \omega'+\ddb u_0+\mu\;\tau^*\omega_1
\in\Omega^{1,1}(\oW).
$$

Recall that the section $s$ defining~$D$ is the pull-back of a
section $s_1$ over~$\CP^1$. We can choose a Hermitian bundle metric
over $\CP^1$ depending on a real parameter~$a$, so that
$e^{a/2}|s_1(z)|=|z|$, for small $z$, and the $(1,1)$-form   
$\Phi=\pbp(\log|s_1|^2)^2$ on
$\CP^1\minus\{z=0\}$ is positive on a punctured neighbourhood of $z=0$,
expressed there as  $2(dz/z)\wedge (d\bar{z}/\bar{z})$.
The ambiguity to choose the constant~$a$ will be needed later.
The pull-back via $\tau$ defines a Hermitian bundle metric $|\cdot|_a$
on $K^{-1}_{\oW}$, so that $|s(y)|_a=|s_1\circ\tau(y)|$, $y\in\oW$.
Respectively, the $(1,1)$-form $\tau^*\Phi$ on $W=\oW\minus D$
is semi-positive near $D$ and positive in the directions transverse to
the fibres of~$\tau$ near~$D$.

With the above arrangements, define
	\begin{equation}\label{a-form}
\omc = \cpt(\mu)+\frac{i}4\tau^*\Phi 
= \omega'+\ddb u_0+\mu\;\tau^*\omega_1+\frac{i}4\pbp(\log|s|_a^2)^2.
        \end{equation}
The last two terms in the right-hand side of~\eqref{a-form} are
pulled back from $\CP^1$ and it is not difficult to see
that $\omc$ is a well-defined positive non-degenerate
$(1,1)$-form, hence a K\"ahler form on~$W$, for any sufficiently large
constant $\mu>0$. We note, for later 
use, that the section $s$ can also be measured in the Hermitian bundle
metric on $K_{\oW}^{-1}$ induced by~$\omega'$. Near~$D$ this latter
metric satisfies $|s(y)|^2=e^{-\rho(y)}|s(y)|_a^2=e^{-\rho(y)-a}|z|^2$,
where $\rho$ is a smooth function of~$y$.
Complementing $z$ to a system of local holomorphic coordinates
near a point in~$D$, we obtain a local expression
$\rho(x)+a=\log\det g'_{i\bar{j}}(x)$, where $g'_{i\bar{j}}$ are the
local coefficients of~$\omega'$. 

To deduce the asymptotic expression for $\omc$,
put $\kappa_I = (\omega' + \tddb u_0)|_D$ and extend $\kappa_I$
to a closed real 2-form on~$U$ via the pull-back of the projection
$U\to D$ along~$z$, determined by the product decomposition
$U\simeq\{|z|<1\}\times D$. The inclusion $D\subset U$ is a homotopy
equivalence, in particular, any 2-cycle in~$U$ is in the homology
class of some 2-cycle in~$D$. The class
$[(\omega' + \tddb u_0)|_U-\kappa_I]\in H^2(U,\RE)$ then is trivial
since it vanishes on every 2-cycle in~$D$. Thus the restriction 
of the form $(\omega' + \ddb u_0)$ to $U$ may be written as
$\kappa_I + d\nu$, for a smooth real 1-form $\nu$ on~$U$ with
$(d\nu)|_D=0$. We may use the product decomposition on $U$ to expand in
powers of~$z$, writing
	\begin{equation*}
\nu=\nu_0+f_0dz+\bar{f}_0d\bar{z}+O(|z|),
	\end{equation*}
for some 1-form $\nu_0$ and function $f_0$ on~$U$
pulled back from~$D$, i.e.\ independent of~$z$.  The form $\nu$ is
determined only up to addition of an exact 1-form (remember that
$H^1(U,\RE)=H^1(D,\RE)=0$). The property $(d\nu)|_D=0$ implies
$d\nu_0=0$, so $\nu_0$ is exact and may be assumed zero.
The last two terms in~\eqref{a-form} are pulled back from $\CP^1$ and
depend only on~$z$. Direct calculation with the change of variable
$e^{-t-i\theta}=z$ then shows that `near the divisor~$D$ at infinity'
the $\omc$ has the asymptotic form~\eqref{asymp} \mbox{with $\lambda=1$}
	\begin{equation}\label{a-cyl}
\smash[t]{\omc|_{U\minus D} =
dt\wedge d\theta+\kappa_I+d(e^{-t}\psi_{\text{cyl}})=
\omega_0+d(e^{-t}\psi_{\text{cyl}}),}
	\end{equation}
where $\psi_{\text{cyl}}$ is a smooth 1-form bounded with all
derivatives on the cylindrical end.
	\begin{remark}
In fact, we can show by further calculation that $\omc$ can be chosen
so that the last term in \eqref{a-cyl} has the decay rate $e^{-2t}$,
rather than $e^{-t}$. However we omit these details as they are not
needed for the later arguments.
	\end{remark}

The Ricci forms of the metrics $\omega'$ and $\omc$ are related by
$\Ric(\omc)-i\pbp\log({\omega'}^3/\omc^3)
=\Ric(\omega')$, and we can
write
	\begin{equation}
	\label{ricci}
\smash[t]{\Ric(\omc)= i\pbp \fcy,
\text{ where }
\fcy = -\log\frac{\omc^3}{{\omega'}^3}-\log|s|^2,}
	\end{equation}
Here $|s|$ is taken in the metric induced by $\omega'$ and
we used the standard expression $\Ric(\omega')=-i\pbp\log|s|^2$
for the Ricci form via the curvature form of $K^{-1}_W$.
From~\eqref{a-form} we have 
$\omc^3=|z|^{-2}((i/2)dz\wedge d\bar{z}\wedge\cpt^2+|z|^2\cpt^3)$
on the neighbourhood $U\minus D$ in~$W$ when $|z|\neq 0$ is small.
Using part of the standard argument making a K\"ahler metric~$\cpt$
Euclidean to order~2 at a given point (e.g.~\cite{GH}, \S 0.7), we can
complement the coordinate $z$ near any point of~$D$ to local
holomorphic coordinates $z_0,z_1,z_2$, where $z_0=z$, in such a way
that $(i/2)dz\wedge d\bar{z}\wedge\cpt^2= 
(1+O(|z|^2)) \prod_{m=0}^2(i/2)dz_m\wedge d\bar{z}_m$.
On the other hand, in the same local coordinates ${\omega'}^3=\det
(g'_{i\bar{j}}) \prod_{m=0}^2(i/2)dz_m\wedge d\bar{z}_m$.  Using also
the previous calculation $|s|^2=|z|^2\det(g'_{i\bar{j}})$
near~$D$, we can write $|s|^2\omc^3/{\omega'}^3=1+O(|z|^2)$.
Passing to the cylindrical coordinates $t,\theta$, we see, by
inspection of~\eqref{ricci}, that
the function $\fcy$ and any derivatives of~$\fcy$
are bounded on~$W$ by a constant multiple of $e^{-2t}$. Here $t$ is
understood to be extended from $U$ to a smooth non-negative function
on $W$ by cutting off to zero away from a neighbourhood of~$D$.

In order to invoke the second part of Proposition~\ref{soln} with
$f=\fcy$, it remains to show that we can satisfy the condition
$\int_W(e^{f_{\text{cyl}}}-1)\omc^3=0$.  Note that as $\fcy$ decays
exponentially to zero along the cylindrical end of~$W$ the integral is
convergent. Recall that we have a choice of
the real parameter~$a$ in the construction of $\omc=\omca$
which does not affect the asymptotic properties achieved above.
By~\eqref{ricci} the summand $e^{\fcy}\;\omca^3={\omega'}^3|s|^{-2}$ is
independent of~$a$. On the other hand,
	\begin{equation}\label{param}
	\begin{split}
\pbp(\log|s|_a^2)^2&=\pbp(a-\log|s|_0^2)^2\\
&=\pbp(\log|s|_0^2)^2-2a\pbp\log|s|_0^2=
\pbp(\log|s|_0^2)^2-2a\tau^*\Upsilon,
	\end{split}
	\end{equation}
where $\Upsilon=\pbp\log|s_1|^2$ denotes the curvature form of the
Hermitian line bundle over $\CP^1$ discussed earlier in this section.
Calculating from \eqref{a-form} and~\eqref{param},
	\begin{align*}
\int_W(e^{f_{\text{cyl}}}-1)\;\omca^3
&=\int_W\bigl({\omega'}^3|s|^{-2}-
  (\omco^3-3\omco^2\,\frac{ai}2\,\tau^*\Upsilon)\bigr)\\
&=\int_W\bigl({\omega'}^3|s|^{-2}-
  (\omco^3-\frac{3ai}2\,{\omega'}^2\;\tau^*\Upsilon)\bigr),
	\end{align*}
using Stokes' Theorem on~$\oW$ for the latter
equality, we find that $\int_W(e^f-1)\;\omca^3$ depends linearly
on~$a$. Note that $|s_1(z)|=e^{-a/2}|z|$ near $z=0$ implies that the
form $\tau^*\Upsilon$ is compactly supported in~$W$. Further, the class
$c_1(\oW)=(2\pi i)^{-1}[\tau^*\Upsilon]$ of the anticanonical bundle
of $\oW$ is semi-positive and $\omega'$ is positive non-degenerate and
we deduce that $\int_W {\omega'}^2\,\tau^*\Upsilon\neq 0$. Therefore,
the integral $\int_W(e^{f_{\text{cyl}}}-1)\;\omca^3$ vanishes for some
value of~$a$.

Now, for the constructed (1,1)-form $\omc$ and the function
$f_{\text{cyl}}$, there is a unique decaying solution $u$ for the complex
Monge--Amp\`ere equation, according to Propositions~\ref{soln}
and~\ref{unique}.
	\begin{prop}\label{decay}
The solution $u$ satisfies $|\nabla^k u|<C_{k,\eps}
e^{-(\lambda-\eps)t}$, for any $\eps>0$, where
the constant $C_{k,\eps}$ is independent of~$t$.
Here $\lambda=\min\{1,\sqrt{\lambda_1(D)}\}$ and $\lambda_1(D)$ is the
smallest positive  eigenvalue  of the  Laplacian  on  $C^\infty(D)$ for
the metric $\kappa_I$.
        \end{prop}

Before going to prove Proposition~\ref{decay} we need to recall
from~\cite{LM},\cite{tapsit} some Fredholm theory
for the Laplacian on asymptotically cylindrical Riemannian manifolds
`with boundary at infinity'. (Note the remark on~\cite{tapsit} after
Theorem~\ref{ric-flat}.)

We begin by introducing a framework of weighted Sobolev spaces for
exponentially decaying functions.  As before, let $t\in C^\infty(W)$
coincide along the cylindrical end of~$W$ with the $\RE$-coordinate
$t=-\log|z|$.  Consider a smooth asymptotically cylindrical metric
on~$W$, more precisely, a metric which for the large values of $t$ is
asymptotic with all derivatives to a product metric on the cylindrical
end $U\simeq \RE_{>0}\times S^1\times D$ of~$W$ ($\omc$ is an example
of such metric). As $t$ tends to infinity, the metric on $\{t\}\times
S^1\times D\subset W$ converges to some metric on $S^1\times D$ and we
shall call the Riemannian manifold $S^1\times D$ with this limit
metric the boundary of~$W$ at infinity.
The weighted Sobolev space $e^{-\delta t}L^p_k(W)$
is the space of all functions $e^{-\delta t}v$ such that
$v\in L^p_k(W)$. By definition, the norm of $e^{-\delta t}v$ in
$e^{-\delta t}L^p_k$ is just the $L^p_k$-norm of~$v$. The definition
extends in the usual way to differential forms.
        \begin{prop}\label{b-diff}
(i) Let $\tW$ be a smooth Riemannian manifold with a smooth
asymptotically cylindrical metric and let $q$ be an integer
$0\le q\le \dim \tW$. Let $\delta>0$ be such
that $\delta^2$ is {\em not} an
eigenvalue of the Laplacian acting on
$(\Omega^{q-1}\oplus\Omega^q)(\p_\infty\tW)$ (on
$C^\infty(\p_\infty\tW)$ in the case $q=0$), where $\p_\infty\tW$
denotes the boundary of~$\tW$ at infinity.
Then the Laplacian on the differential $q$-forms
on~$\tW$ defines, for every $p>1$, $k\ge 2$, a bounded Fredholm
operator $e^{-\delta t}L^p_k(\Omega^q(\tW))\to
e^{-\delta t}L^p_{k-2}(\Omega^q(\tW))$. 
In particular, its image in $e^{-\delta t}L^p_{k-2}$ is {\em closed}.

(ii) Let $\tW$ be as in (i) and $\delta>0$.
Then the Laplacian $\Delta$ on the space of functions
$e^{-\delta t}L^p_{k}(W')$ is injective.
Given $\tilde{f}\in e^{-\delta t}L^p_{k-2}(W')$ with
$\int_{\tW}\tilde{f}=0$, the equation $\Delta v=\tilde{f}$ has a
(unique) solution $v\in e^{-\delta_1 t}L^p_{k}(W')$, where
$\delta_1\le\delta$, $\delta_1<\lambda_1$, $\lambda_1>0$, and
$\lambda_1^2$ is the smallest positive eigenvalue of the Laplacian
$\Delta_{\p}$ on $\p_\infty\tW$.
In particular, if $0<\delta<\lambda_1$ then the image 
$\Delta(e^{-\delta t}L^p_{k}(W'))$ is a codimension-one subspace in
$e^{-\delta t}L^p_{k-2}(W')$ of functions whose integral over $W'$ is
zero.
        \end{prop}
	\begin{pf}
The part (i) is a direct application of~\cite{LM} or~\cite{tapsit}.

For the part (ii), we find by application of Theorem~7.4 of~\cite{LM}
that the index of the Laplacian on $e^{-\delta t}L^p_{k-2}(W')$ (for
any $p>1$, $k\ge 2$) is $-1$ whenever $0<\delta<\lambda_1$. On the
other hand, for any $\delta >0$ the integration by parts is
valid and proves that then any $e^{-\delta t}L^p_k$ function in the
kernel of the Laplacian on~$W$ must be constant
(cf.~\cite{tapsit} p.224).
As there are no non-zero constant functions in $e^{-\delta t}L^p_k(\tW)$
for $\delta >0$, this gives the injectivity.
The integration by parts also shows that the image
$\Delta(e^{-\delta t}L^p_k(\tW))$ consists of functions with zero
integral over $\tW$ and the rest of (ii) then follows.
	\end{pf}
	\begin{remark}
For the values $\delta>\lambda_1$ of the weight parameter, the
index of the Laplacian will be strictly less than $-1$, by direct
application of the change of index formula in \cite{LM}
or~\cite{tapsit}. Then the Laplacian will no longer be onto
the space $\{\tilde{f}\in e^{-\delta t}L^p_{k-2}(W'):
\int_{\tW}\tilde{f}=0\}$ and so the decay estimate for $v$ in
Proposition~\ref{b-diff}(ii) cannot in general be improved.
	\end{remark}
        \begin{pf}[Proof of Proposition~\ref{decay}]
The complete K\"ahler metric $\omc$ on $W$, being an asymptotically
cylindrical metric, has bounded curvature and injectivity radius
bounded away from zero.  Therefore, we
still have Sobolev embedding $L^p_k(W)\subset C^r(W)$ for
\mbox{$r<k-6/p$} (see \cite{aubin}~\S2.7).
We shall be done if for some fixed~$p>1$ we show that
$u\in e^{-(\lambda-\eps)t}L^p_k(W)$, for any $\eps>0$ and
any integer $k$.

Choose $p>1$, such that $2-6/p>1$, and let $\lambda-\eps>0$ be as in
the hypothesis. The boundary of $W$ at infinity is the Riemannian
product $D\times S^1$, where the metric on $D$ is~$\kappa_I$ and $S^1$
is the unit circle.
The complete set of eigenfunctions of the Laplacian for the product
metric is generated as products of the eigenfunctions on the factors
(\cite{BGM}, p.144), so the eigenvalues of the Laplacian on
$C^\infty(D\times S^1)$ are the sums $\lambda_j(D)+n^2$ of those on $D$
and on~$S^1$.
Thus $(\lambda-\eps)^2$ is less than the first positive eigenvalue
of the Laplacian of the product metric on $D\times S^1$.
With the above choice of weight $\lambda-\eps$ and Sobolev parameters
$p$ and $k\ge 2$, the integration by parts can be applied, so that
there is a well-defined non-linear map
        \begin{equation}\label{lapl}
	\begin{split}
\eus{A}: e^{-(\lambda-\eps)t}L^p_{k}(W) &\to
\{\tilde{f}\in e^{-(\lambda-\eps)t}L^p_{k-2}(W):
\int_W\tilde{f}\;\omc^3=0\},
\\
\eus{A}(v)&= 
\frac{(\omc+\tddb v)^3}{\omc^3}-1,
	\end{split}
        \end{equation}
which we further restrict to the open subset of those $v$ satisfying
$\omc+\tddb v>0$, so that \mbox{$\omc+\tddb v$} is a metric.
The map~$\eus{A}$ is smooth with the derivative at $v$ given by 
$(d\eus{A})_v=(3(\omc+\tddb v)^3/(2\pi\omc^3))\Delta_v$ where
$\Delta_v$ is the Laplacian associated to the metric $\omc+\tddb v$
(cf.~\cite{aubin} \S 7.4).  If $v$ is in $e^{-(\lambda-\eps)t}L^p_k(W)$
for any~$k$, with some fixed $p$ and $\eps$, then $\omc+\tddb v$
defines a smooth asymptotically cylindrical metric on~$W$.
In particular, the coefficient before $\Delta_v$ in the derivative
of~$\eus{A}$ is smooth and asymptotically~1 along the cylinder.
Also, by Proposition~\ref{b-diff}, the Laplacian $\Delta_v$ is injective
and maps onto the target space of~$\eus{A}$ in~\eqref{lapl} which is a
closed subspace of the weighted Sobolev space
$e^{-(\lambda-\eps)t}L^p_{k-2}(W)$.
Thus $(d\eus{A})_v$ is an isomorphism between Banach spaces
$e^{-(\lambda-\epsilon)t}L^p_k(W)$ and
$\{\tilde{f}\in e^{-(\lambda-\epsilon)t}L^p_{k-2}(W):
\int_W\tilde{f}\;\omc^3=0\}$.

We calculated that $\Ric(\omc)=i\p\bar{\p}f$, where the
function $f=f_{\text{cyl}}$ is in $e^{-(\lambda-\epsilon)t}L^p_{k-2}(W)$ for
any $p>1$, $k\ge 2$, \mbox{$\epsilon>0$}, and $\int_W(e^f-1)\,\omc^3=0$.
If we could assume the
norm of~$f$ in $e^{-(\lambda-\epsilon)t}L^p_{k-2}(W)$ to be as small as we
like then the Inverse Function Theorem in Banach spaces would give us
a well-defined $\eus{A}^{-1}(e^f-1)\in e^{-(\lambda-\epsilon)t}L^p_{k}(W)$.
Then the unique decaying solution of the complex Monge--Amp\`ere equation
given by Propositions~\ref{soln} and~\ref{unique} would be given by
$u=\eus{A}^{-1}(e^f-1)$, hence the required decay estimate.

Lacking the assumption that $f$ is small, we start from the assertions
of Proposition~\ref{soln} that the solution $u$ of $\eus{A}(u)=e^f-1$
decays to zero along the cylindrical end of~$W$ and has bounded
derivatives of any order. 
We next show that in fact all the derivatives of $u$ decay to zero
along the end of~$W$.  A suitable method may be adapted from the
regularity arguments in \cite{joyce-book}~p.194, here is an outline.
Consider the linear operator $\eus{D}$ defined by 
$$
\eus{D}v=i\pbp v\; \frac{\omega_u^2+\omega_u\omc+\omc^2}{\omc^3},
$$
where $\omega_u= \omc+\tddb u$.
Let $B(y,2R)\subset W$ denote the geodesic ball about $y\in W$ of
fixed, sufficiently small radius~$2R$, in the metric $\omc$. A
positive $2R$ smaller then the injectivity radius of $\omc$ can be
chosen independent of $y\in W$. Then $\eus{D}$ induces, by restriction
and with the help of the exponential mapping at $y$, an elliptic
operator, $\eus{D}_y$ say, on an open domain in~$\CX^3$. Moreover, the
lower bound on the eigenvalues of the matrix of coefficients of the
principal symbol of $\eus{D}_y$ and the upper bounds on the
coefficients of $\eus{D}_y$ on the domain in $\CX^3$ may be taken
independent of~$y$.  Applying a standard elliptic (Schauder) estimate
on an open domain (e.g.~\cite{GT}, Theorems~6.2 and 6.17), putting
$v=u$ and using
$(2\pi)^{-1}\eus{D}u=\eus{A}(u)=e^f-1$, we can show that
	\begin{equation}\label{elliptic}
	\begin{split}
\|u\|_{C^{2+k}(B(y,R))}&< C (\|\eus{D}u\|_{C^k(B(y,2R))}
+\|u\|_{C^0(B(y,2R))})\\
&< C' (\|f\|_{C^k(B(y,2R))}+\|u\|_{C^0(B(y,2R))}).
	\end{split}
	\end{equation}
with the constants $C,C'$ depending only on the metric $\omc$ and the
norm of $u$ in~$C^{2+k}(W)$. The right-hand side of~\eqref{elliptic},
as a function of $y\in W$, decays to zero along the end of~$W$, so
the derivatives of $u$ decay as well. In particular, the metric
$\omega_u$ is asymptotic, with all derivatives, to~$\omc$ along the
end of $W$.

Now the operator $\eus{D}$ is asymptotic to
$3\Delta_0$, the Laplacian of~$\omc$ and, by the theory of~\cite{LM}
Sec.6, $\eus{D}$ defines a Fredholm operator between the same
$e^{-(\lambda-\epsilon)t}$-weighted Sobolev spaces as
$\Delta_0$. Using the identity 
$(i\pbp v)\omc^2=(\Delta_0 v)\omc^3$ and similar ones for the
metrics $\omega_u$ and $\omc+\omega_u$, we calculate
$2(\eus{D}v)\omc^3=(\Delta_0 v)\omc^3+(\Delta_u v)\omega_u^3+
(\Delta_{0,u}v)(\omc+\omega_u)^3$. Then
an integration by parts argument shows that $\eus{D}$ is
$L^2$ self-adjoint (in the metric $\omc$), injective,
has index $-1$ and the same cokernel as~$\Delta_0$. As
$\eus{D}u=2\pi(e^f-1)\in e^{-(\lambda-\epsilon)t}L^p_{k-2}(W)$ we have
$u\in e^{-(\lambda-\epsilon)t}L^p_{k}(W)$, for all~$k$.
        \end{pf}

This completes the proof of~\eqref{asymp}.
\enlargethispage{4.5ex}
\smallskip

Regarding~\eqref{h.vol}, note first that if $\Omega_g$ is a
holomorphic volume form for $\omega_g$ then $z\Omega_g$ extends to a
holomorphic $(3,0)$-form on $U\subset\oW$.
On the other hand, \ $\Omega_0=(-dz/z)\wedge\;(\kappa_J+i\kappa_K)$, 
so $z\Omega_0$ extends to a holomorphic $(3,0)$-form in a different,
product complex structure on the same neighbourhood
$\{|z|<1\}\times D\simeq U$.
As the metric $\omega_g$ is
asymptotic to $\omega_0$, we may without loss assume
that $z\Omega_g=dz\wedge(\kappa_J+i\kappa_K)$ at any point in
$\{0\}\times D$.
The \mbox{3-forms} $\Omega_g$ and $\Omega_0$ are smooth on~$U$ and
straightforward calculation in coordinates shows that
$\Omega_g-\Omega_0$ and all derivatives are $O(e^{-t})$, hence
certainly $O(e^{-(\lambda-\eps)t})$ on the cylindrical end.
We claim that $\Omega_g-\Omega_0=d\Psi$ for some 2-form $\Psi$
on~$\RE_{>0}\times S^1\times D\simeq U\minus D\subset W$.
Indeed, for any 2-cycle $[C]\in H_2(D)$ the integral
$\int_{\{t\}\times S^1\times C}\Omega_g$ is independent of $t$ by
Stokes' Theorem, as $\Omega_g$ is closed. Letting $t$ tend to infinity
we find $\int_{\{t\}\times S^1\times C}(\Omega_g-\Omega_0)=0$,
and since $C$ was any 2-cycle we have $[\Omega_g-\Omega_0]=0$
in $H^3(U\minus D)$, so $\Omega_g-\Omega_0=d\Psi$.
For the decay properties of~$\Psi$, extend $\Psi$ smoothly from
$U\minus D$ to all of $W$,
so that $d\Psi$ is in $e^{-(\lambda-\eps)t}L^p_k\Omega^3(W)$
for any $k$ and any $\eps>0$.
We choose some $\lambda-\eps>0$ and may now apply to $W$
Proposition~6.13 of~\cite{tapsit} 
(cf. also the argument of Lemma~6.11 therein)
to claim that the zero class $d\Psi$ of the `usual' ($C^\infty$) de
Rham cohomology $H^3_{\text{dR}}(W)$ may be represented as the zero
cohomology class of the $e^{-(\lambda-\eps)t}$-weighted Sobolev
completion of the de Rham
complex. Thus $\Psi$ may be taken to be in
$e^{-(\lambda-\eps)t}L^p_{k+1}\Omega^2(W)$ (moreover, in
$e^{-t}L^p_{k+1}\Omega^2(W)$), for every $k$, and the
required decay property follows by Sobolev embedding.

\section{The compact 7-manifolds}
\label{K3}\label{7mfd}

The Riemannian products of $S^1$ and the asymptotically cylindrical
$SU(3)$-manifolds of the previous section give real 7-dimensional
manifolds with holonomy $SU(3)$. These asymptotically cylindrical
7-manifolds have torsion-free $G_2$-structures, by the inclusion
$SU(3)\subset G_2$, and can be `approximated' (in a sense that we
shall make precise below) by compact 7-manifolds with boundary. The
boundary is the Riemannian product of $S^1\times S^1$ and an
$SU(2)$-manifold underlying a Ricci-flat K\"ahler K3 surface. We want
to take the union of two such 7-manifolds, with boundaries identified
via an isometry, to produce a closed compact 7-manifold carrying a
$G_2$-structure with arbitrary small torsion.

A little thought shows that the naive identification of the
boundaries is not a satisfactory option as the closed 7-manifold will
then have an infinite fundamental group because of the $S^1$~factor,
thereby violating a necessary condition for the existence of metrics
of holonomy~$G_2$. We can avoid this topological obstruction by using
a certain `twisted' identification map. The construction of this map
and of the family of `approximate solutions' to the
equations~\eqref{torsion} is explained in this section.

We begin by defining a relation (and an isometry arising from that
relation) for pairs of K3 surfaces. This uses some standard results on
K\"ahler geometry of K3 surfaces which are only briefly recalled here,
for a comprehensive reference see~\cite{BPV}, Ch.VII.

Recall that the $SU(2)$-metric giving the cylindrical asymptotic model
in Theorem~\ref{ric-flat} is determined by a Ricci-flat K\"ahler K3
surface~$D$. As we pointed out earlier, a Ricci-flat K\"ahler metric
on a complex surface is hyper-K\"ahler; in terms of the holonomy
groups this result corresponds to the low-dimensional coincidence
$SU(2)=\SP(1)$. Now a hyper-K\"ahler metric is characterized by the
property that apart from the complex structure $I$ inherited as a
divisor in $\oW$, $D$ admits another two complex structures $J$ and $K$
satisfying the quaternionic relations $IJ=-JI=K$. The $SU(2)$-metric is
K\"ahler with respect to each of the $I,J,K$ (and indeed with respect
to any complex structure in the sphere $aI+bJ+cK$, $a^2+b^2+c^2=1$).
The corresponding K\"ahler forms $\kappa_I,\kappa_J,\kappa_K$ are
orthogonal to each other at every point of $D$ and the volume form of
the hyper-K\"ahler metric may be written as any of the
$\kappa_I^2=\kappa_J^2=\kappa_K^2$. The form $\kappa_J+i\kappa_K$ is a
holomorphic $(2,0)$-form with respect to $I$. Note that there is a
complete $SO(3)$-symmetry between $I,J,K$. In particular, given
a Ricci-flat K\"ahler metric $\kappa_I$ on~$D$, the above normalizations
determine $\kappa_J+i\kappa_K$ only up to a factor $e^{i\vartheta}$,
$\vartheta\in\RE$, i.e.\ there is an $S^1$-ambiguity to choose the
complex structure~$J$ and hence~$\kappa_J$.

Denote by $D_J$ the complex K3 surface defined by considering on
$D$ the complex structure $J$ instead of~$I$. Notice that $(D,\kappa_I)$
and $(D_J,\kappa_J)$ are isometric as the {\em real} 4-dimensional
Riemannian manifolds but are in general {\em not} isomorphic as
complex surfaces.
	\begin{defi}
Let $(D,\kappa_I)$ and
$(D',\kappa'_I)$ be two Ricci-flat K\"ahler
K3 surfaces. We say that these satisfy the {\em matching condition} if
there is a choice of holomorphic $(2,0)$-forms $\kappa_J+i\kappa_K$,
$\kappa'_J+i\kappa'_K$, and an isomorphism of $\ZE$-modules
$h:H^2(D',\ZE)\to H^2(D,\ZE)$
preserving the cup-product and such that $h([\kappa'_I])=[\kappa_J]$,
$h([\kappa'_J])=[\kappa_I]$, $h([\kappa'_K])=-[\kappa_K]$,
for the $\RE$-linear extension of $h$.
	\end{defi}
	\begin{prop}
        \label{torelli}
	\label{match}
Suppose that $(D,\kappa_I,\kappa_J+i\kappa_K)$ and
$(D',\kappa'_I,\kappa'_J+i\kappa'_K)$ satisfy the matching condition.
Then there is an isomorphism of complex surfaces $f:D_J\to D'$, such
that $f^*=h$. Moreover, $f$ is an
isometry of hyper-K\"ahler manifolds, with the pull-back action
on the K\"ahler forms given by
	\begin{equation}
	\label{K3iso}
f^*:\;
\kappa'_I\mapsto \kappa_J,\quad
\kappa'_J\mapsto \kappa_I,\quad
\kappa'_K\mapsto -\kappa_K.
	\end{equation}
	\end{prop}
	\begin{pf}
The first claim is a simple application of the (global) Torelli
theorem for K\"ahler K3 surfaces.
We have $H^{2,0}(D')=\CX[\kappa'_J+i\kappa'_K]$,
$H^{0,2}(D')=\CX[\kappa'_J-i\kappa'_K]$ in the complex structure
of~$D'$, and $H^{2,0}(D_J)=\CX[\kappa_I-i\kappa_K]$,
$H^{0,2}(D_J)=\CX[\kappa_I+i\kappa_K]$ in the complex structure
of~$D_J$. Hence the $\CX$-linear extension $h_\CX$ of the isometry~$h$
given by the matching condition converts Hodge decomposition
$H^2=H^{2,0}\oplus H^{1,1}\oplus H^{0,2}$ for $D'$ into the one
for~$D_J$. As $h$ maps K\"ahler class $[\kappa'_I]$ to a K\"ahler
class $[\kappa_J]$ it satisfies the `effective Hodge isometry'
conditions and so, by \cite{BPV},~Theorem~VIII.11.1, $h$~arises as
a pull-back $h=f^*$ of a uniquely determined biholomorphic map
$f:D_J\to D'$.

The second claim follows as a Ricci-flat K\"ahler metric on a K3
surface is uniquely determined by its K\"ahler class. So the image
$f^*(\kappa'_I)$ of the Ricci-flat K\"ahler metric on $\kappa'_I$
has to be the Ricci-flat $\kappa_J$ on $D_J$.
	\end{pf}

\begin{subequations}\label{su3}
Let $W$ be a complex 3-fold with an $SU(3)$-structure defined by a
K\"ahler form $\omega$ and a holomorphic volume form $\Omega$.
The `product' $G_2$-structure on the real 7-manifold $W\times S^1$ is
expressed by the 3-form $\phi$ in $\Omega^3_+(W\times S^1)$ or its
Hodge dual in $\Omega^4_+(W\times S^1)$, respectively,
	\begin{align}
\phi\, &= \omega\wedge d\theta'
  + \im\Omega,	\label{su3g2}
\\
*\phi &= \pol\omega\wedge\omega
  - \re\Omega\wedge d\theta',	\label{su3g2B}
	\end{align}
where $\theta'$ is a standard parallel 1-form on $S^1$.
The $G_2$-structure $\phi$ satisfies~\eqref{torsion} and defines
on $W\times S^1$ a metric whose holonomy is isomorphic to the holonomy
$SU(3)$ of~$W$. Suppose further that Theorem~\ref{ric-flat} holds
for~$W$ with a product decomposition
$\iota:\RE_{>0}\times S^1\times D\to W$ for the cylindrical end 
of~$W$.
\end{subequations}

Recall from~\S\ref{synopsis}, that if $\phi$ is a $G_2$ structure then
$\phi+\chi$ defines another $G_2$-structure whenever
$\|\chi\|_{C^0}<\epsilon$, in the metric defined by~$\phi$, where
$\epsilon>0$ is an `absolute' constant.
Let $\alpha:\RE\to[0,1]$ denote a cut-off function,
$\alpha(t)\equiv 0$ for $t\le 0$ and $\alpha(t)\equiv 1$ for $t\ge
1$. Using the expressions~\eqref{atinfty} for the
asymptotically cylindrical $SU(3)$-structure on $W$, put
	\begin{align*}
&\omega_T = \omega - d(\alpha(t-T+1)\psi),\\
&\Omega_T = \Omega - d(\alpha(t-T+1)\Psi).
	\end{align*}
Then, as $\psi$ and $\Psi$ decay to zero,
	\begin{equation}\label{cut-off}
\phi_T=\omega_T\wedge d\theta'+\im\Omega_T,
	\end{equation}
gives a well-defined $G_2$-structure
$\phi_T\in\Omega^3_+(S^1\times W)$ for every 
$T\ge T_0>1$, for large enough~$T_0$. 

By the construction, $d\phi_T=0$. But, as the cutting-off construction
of $G_2$-structure changes the Hodge star of the induced metric, the
form $\Theta(\phi_T)$ is in general not closed and the
$G_2$-structure $\phi_T$ has torsion.
It is convenient to measure the torsion of $\phi_T$ using the 4-form
	\begin{equation}
	\label{tau}
\T(\phi_T) \overset{\text{def.}}=
  \Theta(\phi_T) -
  (\pol\omega_T\wedge\omega_T-\re\Omega_T\wedge d\theta'),
	\end{equation}
so $d\T(\phi_T) = d\Theta(\phi_T)$. The form
$(\pol\omega_T\wedge\omega_T-\re\Omega_T\wedge d\theta')$,
is {\em not} the Hodge dual of
$\phi_T$ in the metric defined by $\phi_T$ (it would be the Hodge dual
in the product metric on $W\times S^1$, if one could assume that
$\omega_T$ is a well-defined K\"ahler form of a metric on~$W$).
On the other hand, by comparing the definition of $\phi_T$ with~\eqref{su3}
and~\eqref{su2su3} we see that $\T(\phi_T)$ is supported in the
region $\iota(D\times S^1\times S^1\times[T,T+1])$ where the
cutting-off was performed. Moreover, we have.
	\begin{lemma}
	\label{W}
For any $\eps>0$,
$$
\|\T(\phi_T)\|_{L^p_k(W\times S^1)}<C_{p,k,\eps}\, e^{-(\lambda-\eps)T},
$$
with a constant $C_{p,k,\eps}$ independent of~$T$, where
$p>1$, $k=0,1,2,\ldots$, the exponent $\lambda>0$ for a given $W$ is
as defined in Theorem~\ref{ric-flat}, and the norms are taken in the
metric on~$W\times S^1$ induced by $\phi_T$.
	\end{lemma}
	\begin{pf}
This is deduced in a straightforward way from Theorem~\ref{ric-flat}.
Note that all the differential forms, metrics, and their Hodge stars, 
appearing in the definition of $\T(\phi_T)$ decay like
$O(e^{-(\lambda-\eps)T})$ to the cylindrical model determined by
$\omega_0,\Omega_0$ of~\eqref{su2su3} and the cut-off function is
fixed.
	\end{pf}

Let now $W_1$,$W_2$ be two asymptotically cylindrical manifolds, each
with an $SU(3)$-structure satisfying the assertions of 
Theorems \ref{ric-flat} and~\ref{line}. For $i=1,2$, denote 
$$
M_{i,T} =
(W_i\minus \iota(D_i\times S^1\times ]T+1,\infty[))\times S^1,
$$
so $M_{i,T}$ is a compact Riemannian 7-manifold with boundary
$\p M_{i,T}$ diffeomorphic to $D_i\times S^1\times S^1$.  Suppose that
Proposition~\ref{match} holds for the K3 surfaces $D=D_1$ and
$D'=D_2$. Put on each $M_{i,T}$ the $G_2$-structure $\phi_{i,T}$
as in~\eqref{cut-off}. Define a closed smooth \mbox{7-manifold} $M$ as
a union of $M_{1,T}$ and $M_{2,T}$, with collar neighbourhoods of the
boundaries identified using the diffeomorphism
	\begin{equation}\label{bdry}
	\begin{split}
D_1\times S^1\times S^1\times]T,T+1[&\to
D_2\times S^1\times S^1\times ]T,T+1[,
\\
(y,\theta_1,\theta_2,T+t)&\mapsto (f(y),\theta_2,\theta_1,T+1-t)
	\end{split}
	\end{equation}
where $f$ is the hyper-K\"ahler isometry of K3 surfaces which we
obtained in Proposition~\ref{torelli}.
The $G_2$-structure near the boundary of $M_{i,T}$
is equivalent to the `product' torsion-free \linebreak
{$G_2$-structure} on the cylinder \mbox{$D\times S^1\times
S^1\times\RE$} arising from the $SU(2)$-structure on~$D$.
The corresponding 3-form $\phi_{(D)}$ is expressed, using
\eqref{su2su3} and~\eqref{su3g2}, as
	\begin{equation}\label{neck}
\phi_{(D)} = \kappa_I\wedge d\theta_1 + \kappa_J\wedge d\theta_2
 + \kappa_K\wedge dt + d\theta_1\wedge d\theta_2\wedge dt.
	\end{equation}
The property~\eqref{K3iso} of $f^*$ ensures that $\phi_{(D)}$ is
preserved by the diffeomorphism~\eqref{bdry} and thus the compact
manifold $M$ has a well-defined orientation and a family $\phi_T$,
$T>T_0$, of $G_2$-structures induced from those on $M_1$ and~$M_2$.
(Here we abused notation and started to write $\phi_T$ for a
$G_2$-form on~$M$, but this should not cause confusion.)

The diffeomorphism class of $M$ is independent of~$T$.
We shall write $M_T$ as a short-hand for $M$ endowed with the 3-form
$\phi_T$ (and the induced metric).  One can think of the Riemannian
manifold $M_T$ as a generalized connected sum of $W_1\times S^1$ and
$W_2\times S^2$ with the neck of length approximately~$2T$. The
following theorem summarizes the properties of~$M$ and $\phi_T$ that
we achieved in the above construction.
	\begin{thrm}
	\label{approx}
(i)
One has $b_1(M)=b_1(W_1)+b_1(W_2)$. Further,
$\pi_1(M)=\pi_1(W_1)\times\pi_1(W_2)$; in particular,
$M$ is simply-connected if $W_1$ and $W_2$ are so.

(ii)
There exists $T_0$ such that $\phi_T\in\Omega^3_+(M)$ for $T>T_0$.
Then for any $\eps>0$,
	\begin{equation}
	\label{rhs}
\|\T(\phi_T)\|_{L^p_k(M_T)} < 2C_{p,k,\eps} e^{-(\lambda-\eps)T}
	\end{equation}
in the metric defined by~$\phi_T$, where $C_{p,k,\eps}$ is the constant
appearing in Lemma~\ref{W}. Here $0<\lambda\le 1$ is the smallest of the two
exponential decay parameters (as in Theorem~\ref{ric-flat}), for $W_1$
and~$W_2$.
	\end{thrm}
	\begin{remark}
The topological claims on~$M$ follow by, respectively,
Mayer--Vietoris exact sequence and van Kampen theorem,
applied to $M=M_{1,T}\cup M_{2,T}$.
	\end{remark}

\section{The gluing theorem}
\label{deform}

When the fundamental group of a closed 7-manifold~$M$ is finite, any
torsion-free $G_2$-structure, i.e. a solution~$\phi$ to the
system of equations  
	\begin{equation}
	\label{pdes}
d\phi=0,\qquad d\Theta(\phi)=0,
	\end{equation}
will make $M$ into a Riemannian manifold with holonomy~$G_2$,
\cite{joyce-book} pp.244--245.
We have constructed, in the previous section, a family
$\phi_T\in\Omega^3_+(M)$ of
$G_2$-structure 3-forms satisfying the first of the
equations~\eqref{pdes} and having arbitrary small $d\Theta(\phi_T)$ as
$T$ tends to infinity.  This sets the scene for application of the
perturbative methods to find solutions to~\eqref{pdes}.
It is convenient to look for the solutions in the form
$\phi_T+d\eta$, for an unknown \mbox{2-form} $\eta$.

The equations~\eqref{pdes} are invariant under the action of the group
of diffeomorphisms of~$M$ whose Lie algebra is given by the vector
fields on~$M$. Respectively, the infinitesimal action on a closed form
$\phi_T$ is just the Lie derivative $d(X\lrcorner\phi_T)$ in the
direction of a vector field $X$. So {\it a priori} a solution
$\eta$ would be transverse to the subspace of all the 2-forms
$X\lrcorner\phi_T$ generated by vector fields on~$M_T$.

A way to eliminate the infinite-dimensional symmetry of~\eqref{pdes} and
reduce the torsion-free condition $d\Theta(\phi+d\eta)=0$
to an elliptic equation in~$\eta$ has been worked out
in~\cite{joyce-book}, \S10.3. It uses some special geometry of
$G_2$-structures and the local expression
$\Theta(\phi+d\eta)=\Theta(\phi)+(d\Theta)_\phi d\eta+R(d\eta)$
with the second order remainder $R(d\eta)$ satisfying a quadratic
estimate
	\begin{equation}
	\label{qdr}
\|dR(\chi_1)-dR(\chi_2)\|_{L^p} <
(C_1+C_2\|d\Theta(\phi_T)\|_{L^p_1})
(\|\chi_1\|_{L^p_1}+\|\chi_2\|_{L^p_1})\|\chi_1-\chi_2\|_{L^p_1},
	\end{equation}
for $p>7$ and with the constants $C_1,C_2$ independent of
$\phi_T$ or the manifold.
The result may be stated as follows.
	\begin{prop}[cf.~\cite{joyce-book}, Theorem 10.3.7]\label{prelim}
Let $M'$ be a compact Riemannian 7-manifold whose metric comes from a
$G_2$-structure defined by a smooth closed form $\phi'\in\Omega^3_+(M)$.
Denote by $\langle\cdot,\cdot\rangle$ the induced inner product on
the fibres of vector bundles associated to~$TM'$.

There exist constants $\rho>0$, $\tilde{\rho}>0$ independent
of $M'$ or $\phi'$ and such that if a 4-form $\T$ on~$M'$,
satisfies $d\T=d\Theta(\phi')$ and $\|\T\|_{C^0}<\tilde{\rho}$,
then
for any 2-form $\eta$ on~$M'$, with $\|\eta\|_{C^1}\le\rho$,
the equation
	\begin{equation}
	\label{laplace}
(dd^*+d^*d)\eta +
*d\bigl( (1 + {\textstyle\frac13}\langle d\eta,\phi'\rangle)
\T\bigr)  - *dR(d\eta) = 0,
	\end{equation}
implies $d\Theta(\phi'+d\eta)=0$.
	\end{prop}
In view of the estimate~\eqref{rhs} the condition
$\|\T(\phi_T)\|_{C^0}<\tilde{\rho}$
in the above proposition will be satisfied for every sufficiently
large~$T$. Thus we can now turn to solving~\eqref{laplace} for an
unknown small $\eta$ putting $\phi'=\phi_T$ and $\T=\T(\phi_T)$
on the compact 7-manifold~$M'=M_T$.

The terms involving second order derivatives of~$\eta$
in~\eqref{laplace} are the Laplacian term, the linear term in $\eta$ with
coefficients controlled by $\T$, and the quadratic term which may be
regarded as a linear operator whose coefficients depend on $\eta$ and
are bounded by $\sup_{M_T}|d\eta|$. This tells us that~\eqref{laplace}
is an elliptic problem, at least for a $C^1$-small $\eta$ and with $T$
assumed large. If $d\eta$ is H\"older continuous and small in
$C^0(M_T)$ then elliptic regularity arguments show that a
weak, $L^p_2$ solution $\eta$ of~\eqref{laplace} is in fact smooth.
(We do not include the details for this latter claim but see
\cite{joyce-book} p.303, cf. also~\cite{yau} p.363.)
Recall that for $\phi_T+d\eta_T$ to be a well-defined $G_2$-structure,
the term $d\eta_T$ must be small in the sup-norm on~$M_T$. If we are
to solve the second order equation~\eqref{laplace} for $\eta$ in
Sobolev space $L^p_2$, then we need an embedding
$L^p_2\subset C^1$ in dimension~7, which requires $p>7$. In what
follows we fix such a choice of~$p$.

The existence of torsion-free $G_2$-structures on $M$ will be
established by the following technical result.
	\begin{thrm}
	\label{perturb}
There exist
$T_0>0$ such that for each $T\ge T_0$ the equation~\eqref{laplace}
has a unique smooth solution 2-form $\eta_T$ on~$M_T$ satisfying
$\|\eta_T\|_{L^p_2}<K_{p,\delta} e^{-(\lambda-\delta)T}$ and
$\|\eta_T\|_{C^1}<K_\delta e^{-(\lambda-\delta)T}$, for any positive
$\delta$, where $K_{p,\delta},K_\delta$ are constants
independent of~$T$. Here $\lambda$ is a constant determined in
Theorem~\ref{approx}.
	\end{thrm}

Theorem~\ref{perturb} may be viewed as an instance of a gluing theorem
for solutions of non-linear elliptic PDEs.
We should remark right away that the equation~\eqref{laplace} for
torsion-free $G_2$-structures was solved in~\cite{joyce-book},
\S\S11.6-11.8 in a different geometrical setting. In principle, our
gluing theorem could be deduced by verifying that the $G_2$-structures
$\phi_T$ with small torsion can be written so that they formally
satisfy the hypotheses of Theorem~11.6.1 in~\cite{joyce-book} and then
invoke that theorem to claim the existence of solutions
to~\eqref{laplace}.

We shall prove Theorem~\ref{perturb} via a different, more
direct and geometrically natural approach tailored to generalized
connected sums. Another reason for our choice is that this proof
follows up the applications of elliptic theory and Hodge theory on
Riemannian manifolds with asymptotically cylindrical ends, already
required in~\S\ref{bCY}. We also are able to explicitly eliminate
the issue of the `obstruction spaces' (whose vanishing is sufficient
for the existence of a solution) in this gluing problem.

The remainder of this section deals with the proof of Theorem~\ref{perturb}.
The following result is a standard variant of the contraction mapping
principle and is included here for the reader's convenience and to 
keep track of the relations between various constants in the estimates.
	\begin{prop}\label{IMT}
Let $a:\BE \to \BF$ be a smooth map between Banach spaces, $a(\eta) =
a_0 + A\eta + Q(\eta)$. Suppose that $A$ is a linear isomorphism
between $\BE$ and~$\BF$
and that for all $\eta_1,\eta_2\in \BE$
	\begin{equation}
	\label{quadratic}
\| A^{-1}Q(\eta_1) - A^{-1}Q(\eta_2) \|  \leq  
m (\|\eta_1\| + \|\eta_2\|)\|\eta_1 - \eta_2\|  .
	\end{equation}
If
	\begin{equation}
	\label{bounded}
\| A^{-1} a_0 \| < 1/(8m)
	\end{equation}
then the equation $a(\eta) = 0$ has a unique solution~$\eta_0$ of norm
less than~$1/(4m)$.
This solution also satisfies $\|\eta_0\|\le 2\|A^{-1}\|\cdot\|a_0\|$.
	\end{prop}
In our situation, $a_0=-\!*\!d\T(\phi_T)$ is exponentially small
in~$T$, as $T\to\infty$, by Theorem~\ref{approx}(ii). Hence for
large~$T$ the term $C_2\|d\Theta(\phi_T)\|=C_2\|d\T(\phi_T)\|$ in the
quadratic estimate~\eqref{qdr} of $Q(\eta)=\mbox{$-*\!dR(d\eta)$}$ can
be absorbed by taking a larger constant~$C_1$ independent of~$T$.
Furthermore, the linearization $A=A_T$ of~\eqref{laplace} is the
Laplacian on $M_T$ perturbed by adding a linear differential operator
$\textstyle\frac13\! *\! d(\langle d\eta,\phi_T\rangle\T(\phi_T))$
whose coefficients are exponentially small in~$T$. Therefore, for
every sufficiently large~$T$, the operator $A_T$ will have an elliptic
symbol and will be invertible on the orthogonal complement
$(\HH^2)^\bot$ of harmonic 2-forms on~$M_T$. 
As we are really interested in the value of $d\eta$, rather than
$\eta$, and as all the terms in the left-hand side of~\eqref{laplace}
take values in $(\HH^2)^\bot$, we may consider~\eqref{laplace} as an
equation in~$(\HH^2)^\bot$. 
The main issue in the application of Proposition~\ref{IMT}
to~\eqref{laplace} is to give a suitable upper bound on the operator
norm of $A_T^{-1}$ for any large~$T$, so as to satisfy~\eqref{quadratic}
and~\eqref{bounded} together.

Our estimate for the linearized problem on $M_T$ will be
obtained by application of the technique developed in 
\cite{asdsums}.
The basic principle is that the elliptic differential operator $A_T$ 
may be obtained by gluing together two elliptic operators $A_j$ on the
non-compact manifolds $W_j\times S^1$, $j=1,2$ at their asymptotically
cylindrical ends. More precisely, the coefficients of the operator
$A_T$ depend smoothly
on the \mbox{$G_2$-structure} $\phi_T$ on~$M_T$, and the operators
$A_j$ are determined, in just the same way,
by the $G_2$-structures which we used in the construction of $\phi_T$.
As each of the non-compact manifolds $W_j\times S^1$ was given a
torsion-free $G_2$ structure, we may take $\T_j=0$, so the associated
linear operators $A_j$ are just the Laplacians for the respective
metrics. We define
	\begin{align*}
A_1=dd^*+d^*d:&\quad
e^{-\delta t_1}L^p_2\Omega^2(W_1\times S^1)\to
e^{-\delta t_1}L^p\Omega^2(W_1\times S^1),\\
A_2=dd^*+d^*d:&\quad
e^{\delta t_2}L^p_2\Omega^2(W_2\times S^1)\to
e^{\delta t_2}L^p\Omega^2(W_2\times S^1),
	\end{align*}
where $\delta>0$ and, as before, $t_j\in C^\infty(W_j\times S^1)$
denotes a real parameter along the cylindrical end. From
Proposition~\ref{b-diff}(i), each $A_j$ is a bounded Fredholm operator
if the weight parameter $\delta>0$ is small.
We shall further need the following.
	\begin{prop}\label{inj}
For every sufficiently small $\delta>0$, (i) the Laplacian $A_1$ is
injective, and (ii) every non-zero element of the kernel of the
Laplacian $A_2$ is smooth and has an asymptotic expansion as
$t_2\to\infty$ with the leading term polynomial in~$t_2$.
	\end{prop}
	\begin{pf}
This is an application of a general result on the asymptotic
expansions proved in \cite{tapsit}~Proposition~5.61.
The leading term of asymptotic expansion of an
$e^{\epsilon t}L^p_k$ kernel element may be written as
$e^{-\mu t_j}\eta(t_j,y)$, where $-\mu<\epsilon$ and a form
$\eta(t,y)$ depends polynomially on $t_j$ and smoothly on the
coordinates~$y$ of the cross-section of the ends of~$W_j\times S^1$.
Further, $\mu^2$ is an eigenvalue of the Laplacian (on differential
forms) on the cross-section at infinity (pp.224--225 op.cit.).
The cross-sections of $W_j\times S^1$ at infinity are isometric;
choose $\delta>0$ with $\delta^2$ smaller than the first positive
eigenvalue of the respective Laplacian.

For $\epsilon=-\delta$, it follows that $\Ker A_1$ coincides with the
$L^2$-kernel. 
It is proved in \cite{aps},~Proposition~4.9, and
\cite{tapsit},~Proposition~6.14, that the $L^2$-kernel of the
Laplacian on differential forms on an asymptotically cylindrical
manifold is isomorphic to the image of the natural
homomorphism $H^*_c\to H^*$ of the de Rham cohomology of this
manifold, where $c$ indicates the cohomology with compact support.
As $H^1_c(W_j,\RE)$ and $H^1(W_j,\RE)$ vanish, it suffices to consider
the image of $H_c^2(W_j,\RE)\to H^2(W_j,\RE)$.
The latter map is a part of exact sequence 
$H^2_c(W_j,\RE)\to H^2(W_j,\RE)\to H^2(D\times S^1,\RE)$. The map of
$H^2(W_j,\RE)$ is equivalent to the composition of injective homomorphisms
$H^2(W_j,\RE)\to H^2(\oW_j,\RE)\to H^2(D,\RE)$, the first one from the
Mayer--Vietoris sequence, as in Theorem~\ref{line}, the second by
Lefschetz hyperplane theorem. (Note that $H^2(D\times S^1,\RE)\cong
H^2(D,\RE)$.) It follows that $H^2_c(W_j\times S^1,\RE)\to
H^2(W_j\times S^1,\RE)$ maps to zero and $A_1$ is injective.

If $\epsilon=\delta$ then, by the same argument as above, a non-zero
2-form in the kernel of $A_2$ cannot be in $e^{-\delta t_2}L^p_2$.
So, by the choice of~$\delta$, the only possibility for the leading
exponent in the asymptotic expansion of such 2-form is $\mu=0$.
	\end{pf}

We now define the weight function $w_T\in C^\infty(M_T)$. Let
$-T<t<T$ denote the real parameter along the neck of~$M_T$, increasing
towards the $(W_2\times S^1)$-piece,
and with $t=0$ in the middle of the cutting-off region
(recall~\S\ref{7mfd}). The $w_T$ coincides along the neck with
$e^{-\delta t}$ and is smoothly cut off to the constants
$e^{\pm\delta T}$ away from the neck, so that
$e^{-\delta T}\le w_T\le e^{\delta T}$ at every point in~$M_T$. Thus
$w_T$ interpolates between the weights on $W_j\times S^1$.
The weighted $w_T L^p_k$-norms on a compact $M_T$ are commensurate,
for any $T$, to the usual $L^p_k$ Sobolev norms, however the estimates
depend on~$T$,
$e^{-\delta T}\|\eta\|_{p,k}\le\|w_T\,\eta\|_{p,k}
\le e^{\delta T}\|\eta\|_{p,k}$.
Put $\BE_T=(\HH^2)^\bot\cap L^p_2\Omega^2(M_T)$,
$\BF_T=(\HH^2)^\bot\cap L^p\Omega^2(M_T)$,
but use the {\em weighted} norms on these Banach spaces,
$\|\eta\|_{\BE_T}=\|w_T\,\eta\|_{L^p_2}$,
$\|\eta\|_{\BF_T}=\|w_T\,\eta\|_{L^p}$.
Finally, define
	\begin{equation}\label{AT}
A_T:\eta\in \BE_T\to (dd^*+d^*d)\eta + {\textstyle\frac13}
*\!d\bigl( \langle d\eta,\phi_T\rangle \T(\phi_T)
\bigr)  \in \BF_T .
	\end{equation}

	\begin{prop}\label{linest}
Suppose that the weight parameter $\delta>0$ is small, so that the
Laplacians $A_j$ are Fredholm operators and satisfy the assertions of
Proposition~\ref{inj}. Then there exists $T_*>0$ such that for any
$T>T_*$ the inverse of the operator $A_T$ defined by~\eqref{AT} exists
and satisfies $\|A_T^{-1}\|<Ge^{\delta T}$ with a constant $G$
independent of~$T$.
	\end{prop}
	\begin{pf}
We checked above that the $A_T$ is invertible for every
sufficiently large~$T$. It suffices to prove that there exist
constants $G,T_*$, such that $\|A_T\eta\|>Ge^{-\delta T}\|\eta\|$, for
$\eta\in\BE_T$, whenever $T>T_*$. The proof goes by contradiction. Assume
that Proposition~\ref{linest} is {\em not} true, so there exists a sequence
$\eta_n\in\BE_{T_n}$, with $\|\eta_n\|=1$, $T_n\to\infty$, such
that $e^{\delta T_n}\|A_{T_n}\eta_n\|\to 0$.

This is an example of the situation considered in~\cite{asdsums}~\S4.1;
the only difference is that in the present application we do not avoid
the `asymptotic kernels and cokernels' and claim a slightly weaker
estimate.  Let $\beta_n=\beta\circ t_n\in C^\infty(M_T)$, where $t_n$
is the real parameter along the neck of $M_{T_n}$ as above, and
$\beta$ be a cut-off function equal to 0 on $(-\infty,1]$ and to 1 on
$[2,\infty)$. The arguments of Proposition~4.2 and Lemma~4.7
of~\cite{asdsums} yield
	\begin{equation}\label{ks}
\|(1-\beta_n)\eta_n\|\to 0 \quad\text{and}\quad
e^{\delta T_n}\|A_2(\beta_n\eta_n)\|\to 0,
	\end{equation}
Here for the first statement we also took into account that the
operator $A_1$ has no kernel. In the second statement, we identified
$\beta_n\eta_n$ with a form on $W_2\times S^1$.

Let $\beta_n\eta_n = \eta'_n+\eta''_n$, respective to the
$L^2$-orthogonal decomposition $\Ker A_2 \oplus (\Ker A_2)^\bot$ of
the domain of~$A_2$. Then $\eta''_n\to 0$ from~\eqref{ks}. Also,
since the kernel of $A_2$ is finite-dimensional, we may assume,
passing to a  subsequence if necessary, that
$\eta'_n\to \eta_0\in\Ker A_2$.
In view of the asymptotic expansion found in Proposition~\ref{inj},
we have $\|e^{-\delta t_2}\eta_0\|_{L^p_2\{t_2>T\}}\ge C_0e^{-\delta T}$,
with $C_0$ independent of $T>T_*$ and $C_0=0$ if and only if $\eta_0=0$.
Using finite dimension of $\Ker A_2$ again, we can estimate by
standard methods  $\|e^{-\delta t_2}\eta'_n\|_{L^p_2\{t_2>T\}}\ge
C'_0 \|e^{-\delta t_2}\eta_0\|_{L^p_2\{t_2>T\}}$, with $C'_0\neq 0$
independent of~$T>T_*$. It follows that
$$
\|A_2(\beta_n\eta_n)\|=\|A_2(\eta''_n)\|\ge
C\|e^{-\delta t_2}\eta''_n\|_{L^p_2\{t_2>T_n\}}=
C\|e^{-\delta t_2}\eta'_n\|_{L^p_2\{t_2>T_n\}}>
CC'_0C_0 e^{-\delta T_n},
$$
since $\beta_n\eta_n=0$ on $\{t_2>T_n\}$ and where $C\neq 0$. Then
$C_0=0$ from~\eqref{ks}, hence $\eta_0=0$ and $\beta_n\eta_n\to 0$. This
leads to a contradiction, 
$1=\|\eta_n\|\le\|(1-\beta_n)\eta_n\|+\|\beta_n\eta_n\|\to 0$,
which proves Proposition~\ref{linest}.
	\end{pf}

We can replace the $L^p_k$ norms in the quadratic
estimate~\eqref{qdr} by the weighted $w_TL^p_k$ norms, at the cost
of weakening the right-hand side by the factor $e^{-3\delta T}$.
Likewise,\linebreak
$\|d\T(\phi_T)\|_{\BF_T} \le
e^{\delta T}\|d\T(\phi_T)\|_p$ and from Proposition~\ref{linest}
the action of $A_T^{-1}$ amounts to an extra factor $e^{\delta T}$
in the estimates. We find that Proposition~\ref{IMT} applies with
$\BE=\BE_T$, $\BF=\BF_T$, $A=A_T$, and $m = \const\cdot\, e^{4T\delta}$
to solve~\eqref{laplace} for every large $T$, provided that
$0<6\delta<\lambda$, where $\lambda$ is the constant in~\eqref{rhs}.
This completes the proof of Theorem~\ref{perturb}.

\section{From Fano 3-folds to $G_2$-manifolds}
\label{fano}

Let $V$ be a 3-dimensional complex manifold with $c_1(V)>0$. It
follows from Yau's proof of the Calabi conjecture~\cite{yau} that $V$
has a K\"ahler metric of positive Ricci curvature. Consequently, $V$
is simply-connected by a theorem of Kobayashi~\cite{kobayashi} and
$H^2(V,\CX)=H^{1,1}(V)$ by Bochner's theorem. It then
follows that $V$ is projective, by the Kodaira embedding theorem. The
condition $c_1(V)>0$ is also equivalent to the anticanonical sheaf
$\OH(-K_{V})$ being ample, and any such manifold $V$ is called a
\mbox{\em Fano 3-fold}.

For every Fano 3-fold $V$, a generic divisor $D$ in the anticanonical
linear system $|-K_V|$ is a K3 surface, by a theorem of
Shokurov~\cite{sh1}. The self-intersection class of~$D\subset V$ is
Poincar\'e dual to the restriction $c_1(V)|_D$ and therefore cannot be
trivial. However, we have
	\begin{prop}
	\label{blowup}
Suppose that a K\"ahler 3-fold $V$ and a surface $D\in|-K_V|$ are
given. Let $C$ be a smooth curve in $D$ representing $D\cdot D$ and
$\sigma:\tV\to V$ the blow-up of $V$ along~$C$.

Then the closure $\tD$ of $\sigma^{-1}(D\minus C)$ is a smooth
anticanonical divisor on $\tV$ and \mbox{$\tD\cdot\tD=0$.} Further,
$\sigma$ restricts to give an isomorphism $\tD\to D$ of complex
surfaces and this isomorphism identifies the induced K\"ahler metric
on $D\subset V$ with the restriction to $\tD$ of some K\"ahler metric
on~$\tV$.
	\end{prop}
	\begin{pf}
See e.g.~\cite{GH}, pp.608--609. The formula for the anticanonical
class of a blow-up yields, for a curve in a 3-fold,
$|-K_{\tV}|=|\sigma^*(D) - E|=|\tD|$, where $E$ is the exceptional
divisor. It is clear that $\sigma_*$ identifies the cycle
$\tD\cdot\tD$ with a cycle on~$V$. 
But $\sigma_*(\tD\cdot\tD)=0$ (cf. blowing-up a surface at a point).
For the claim on K\"ahler metrics, it can be shown in a standard way,
cf.\ pp.186--187 op.cit., that there is a closed semi-positive
$(1,1)$-form $\omega_{[E]}$ representing $c_1([E])$ with the following
properties: $\omega_{[E]}=0$ on $\tD$ and if $\omega_V$ is a positive
$(1,1)$-form on~$V$ then $\sigma^*\omega_V-k^{-1}\omega_{[E]}$ is
positive on~$\tV$ for any sufficiently large~$k$.
	\end{pf}
        \begin{cor}\label{bFano}
Let $V$ be a Fano 3-fold, $D\in|-K_V|$ a K3 surface, and $\tV$ the
blow-up of~$V$ along a curve $D\cdot D$ as in Proposition~\ref{blowup}.
Then $\tV\minus\tD$ has a complete smooth metric with holonomy $SU(3)$,
making $\tV\minus\tD$ into a manifold with asymptotically
cylindrical end, as defined in Theorem~\ref{ric-flat}. Also,
$\pi_1(\tV\minus\tD)=1$. 
	\end{cor}
	\begin{pf}
An exceptional curve $\ell=\sigma^{-1}(x)$ in $\tV$ (where $x\in C$)
meets $\tD$ in exactly one point, $\ell\cdot\tD=1$. 
The corollary now follows from Theorem~\ref{line} and
Proposition~\ref{blowup}, putting $\oW=\tV$ (recall that for a
blow-up, $\pi_1(\tV)=\pi_1(V)=1$).
	\end{pf}

Consider now two pairs $(\tV_1,\tD_1),(\tV_2,\tD_2)$, where each pair
is obtained by blowing up a curve $D_j\cdot D_j$ in a Fano
3-fold~$V_j$ and lifting a K3 divisor $D_j\in|-K_{V_j}|$ via the
proper transform, as above. Put $W_j=\tV_j\minus\tD_j$. The connected
sum construction in~\S\ref{7mfd} for the pair $W_1\times S^1$ and
$W_2\times S^1$ requires the matching condition to hold for the pair
of K3 surfaces $\tD_j\subset\tV_j$ and their K\"ahler classes.  In
view of Proposition~\ref{blowup} and Corollary~\ref{bFano}, the
possibility to satisfy the matching condition for $\tD_1$ and $\tD_2$
can be decided by considering the divisors $D_1$ and $D_2$ on the Fano
3-folds.  In general, the K\"ahler classes induced on $D_1$ and $D_2$
by the respective embeddings, may force different volumes of $D_1$
and~$D_2$. So it may be necessary to rescale the projective metrics
on $V_1$ or $V_2$ by a constant factor to achieve the matching K3
divisors.

If the two K\"ahler K3 surfaces $D_j$ do match, and so a compact real
7-manifold $M$ and a $G_2$-structure $\phi_T$ on $M$ are well-defined,
then $b_1(M)=0$ by Theorem~\ref{approx}(i). Hence $M$ has a metric
with holonomy $G_2$, by Theorem~\ref{perturb}.  As a short-hand we
shall say that the compact \mbox{$G_2$-manifold}~$M$ is {\em
constructed from the pair of Fano 3-folds} $V_1$ and $V_2$.
This presumes, as the initial step, appropriate choices of the
anticanonical K3 divisors $D_1$ and~$D_2$ and rescaling of the metrics
on $V_1$ and $V_2$ if necessary.

Here is the main theorem of this section.
By a deformation class we mean a maximal family of smooth complex
manifolds which are deformations of each other.
	\begin{thrm}\label{ult}
Let $\curly{V}_1$ and $\curly{V}_2$ be two (not necessarily distinct)
deformation classes of Fano 3-folds. There exist representatives
$V_j\in\curly{V}_j$, $j=1,2$, such that a compact \mbox{$G_2$-manifold}
can be constructed from the pair $V_1$ and $V_2$.
	\end{thrm}
	\begin{remarks}
(i) There is a complete classification of Fano 3-folds into 104
deformation classes \cite{Isk77,mukaimori}, see also Chapter~12
in~\cite{fanovar} for a summary list including the values of
basic invariants. This leads to an interesting question of finding the
`geography' of the respective $G_2$-manifolds, e.g.\  plotting their
topological invariants as in~\cite{joyce-book}, Figure~12.3. The issue
is not taken up here, however we do discuss in \S\ref{examples}
concrete examples of the {\em new} topological types of
$G_2$-manifolds obtained by Theorem~\ref{ult}.

(ii) We also point out that it is not absolutely
necessary to start from quasiprojective varieties obtained from
blow-ups of Fano 3-folds. Indeed, the reader will notice that we
really only needed our 3-folds $\tV$ to be fibred over a Riemann
surface with a K3 fibre in the anticanonical class and that the
complement of a fibre have finite fundamental group. Examples of
$G_2$-manifolds will arise from any pair of such fibred spaces
provided that a K3 fibre in each space is chosen so that the matching
condition is satisfied. It is likely that further examples of
matching pairs of K3 fibrations will also be found, e.g.\
using singular Fano varieties.
	\end{remarks}
As one may expect from the previous arguments, the proof of
Theorem~\ref{ult} is largely a matter of finding the matching
anticanonical K3 divisors. Before dealing with that, we need to identify
the type of K3 surfaces which occur in the anticanonical linear
systems on Fano 3-folds in a given deformation class~$\V$.

If $V$ is a Fano 3-fold then any K3 surface $D$ in $|-K_V|$ is
necessarily projective-algebraic.  Given a projective embedding
of~$V$, it is clear that the degree of an algebraic K3 surface in the
anticanonical linear system is determined by the deformation class of
Fano 3-fold, but in fact more it true.
By Lefschetz hyperplane theorem, in the form due to Bott~\cite{bott},
the embedding $\iota:D\hra V$ induces an injective homomorphism
$\iota^*:H^2(V,\ZE)\to H^2(D,\ZE)$ with $H^2(D,\ZE)/\iota^*H^2(V,\ZE)$
torsion-free. The latter property follows by standard topology as both
$D$ and $V$ are simply-connected and up to a homotopy equivalence $V$
is obtained from $D$ by attaching cells of real dimension at least~3. 
Similarly, the K\"ahler class of $V$ defined by a projective
embedding is primitive (indivisible by an integer $>1$) and the
induced K\"ahler class of~$D$ is primitive
in~$H^2(D,\ZE)$. Recall that a sublattice is called
primitive when the quotient of this sublattice has no torsion. Thus
for every K3 anticanonical divisor $D$ on a Fano 3-fold $V$ in a given
deformation class $\curly{V}$, the Picard lattice $H^2(D,\ZE)\cap
H^{1,1}(D)$ necessarily contains a primitive sublattice
$S(\V)=\iota^*H^2(V,\ZE)$. We use the notation $S(\V)$ to indicate
that this sublattice is determined by the deformation class of~$V$. (We
remind that the K3 surfaces in $|-K_V|$ form a Zariski open, hence
connected, family.)

The algebraic K3 surfaces whose Picard lattice contains a primitively
embedded non-degenerate lattice $S$ are studied
in~\cite{dolgachev} and called `ample $S$-polarized K3 surfaces' if
$S$ contains a K\"ahler class. The cup-product on Picard lattice of
any algebraic surface has positive index exactly 1, by Hodge index
theorem. As a sublattice of the form $S=S(\V)$ contains an ample
class, it is non-degenerate of signature $(1,t)$, where the rank
$1+t=b^2(V)$ is equal to the Betti number of some (any) $V\in\V$. By
the classification result, any Fano 3-fold has $1\le b^2(V)\le 10$ so the
possible values of $t$ are $0\le t\le 9$. Thus the inclusion map
defines an ample $S(\V)$-polarization for every anticanonical K3
divisor on $V\in\V$.

Recall that the integral second cohomology of a K3 surface considered
with the cup-product is isometric to the (unique) even unimodular
lattice $L$ of signature $(3,19)$. In fact,
$L=(-E_8)\oplus(-E_8)\oplus H\oplus H\oplus H$, where
$H=\bigl(\begin{smallmatrix} 0&1\\ 1&0 \end{smallmatrix}\bigr)$
denotes the hyperbolic plane lattice and $E_8$ is a standard (even and
positive-definite) root-lattice. The constraint on the signature of
$S(\V)$ found above implies that a primitive embedding of any $S(\V)$
into~$L$ is unique up to an isometry of~$L$, by
\cite{dolgachev}~Corollary~5.2 (see also \cite{nik}~Theorem~1.14.4).
A map $f:D\to D'$ is defined to be an isomorphism of $S(\V)$-polarized
K3 surfaces if $f$ is an isomorphism of complex K3 surfaces and $f^*$
intertwines the embeddings of $S(\V)$ in the respective Picard
lattices. In combination with the global Torelli theorem this tells us
that the class of $S(\V)$-polarized K3 surfaces is determined by the
deformation class~$\V$ of Fano 3-folds. We shall denote the moduli
space of all the isomorphism classes of $S(\V)$-polarized K3 surfaces
by $\K3(\V)$.

When $b^2(V)=1$, the lattice $S(\V)$ is generated by one element, an
integral (K\"ahler) class on $D$ of positive square say $2n-2$ ($n\ge
2$). In this case $\K3\,(\V)=\K3\,(n)$ is just the familiar
19-dimensional moduli space of algebraic K3 surfaces of degree $2n-2$
in $\CP^n$.  The general case is quite similar. Firstly, by
application of the embedding result of Nikulin (\cite{nik},
Theorem~1.12.4), any lattice $S(\V_1)$ admits a primitive embedding
into the even unimodular lattice
$L_0=(-E_8)\oplus(-E_8)\oplus H\oplus H$ of signature $(2,18)$,
because $S(\V)\otimes\RE$ embeds in $L_0\otimes\RE$ and
$\rank S(\V)\le\pol\rank L_0$. Therefore, the orthogonal complement to
$S(\V_1)$ in the K3 lattice~$L$ contains a copy of~$H$. It then
follows from \cite{dolgachev}~Theorem~5.6
that the moduli space $\K3\,(\V)$ is an irreducible (hence connected)
quasiprojective algebraic variety of (complex) dimension $20-b^2(V)$.

We shall need the following result, which may be viewed as a kind
of converse to Shokurov's theorem cited earlier.
	\begin{thrm}\label{onto}
Let $\curly{V}$ be a deformation class of Fano 3-folds.
Denote by $\MM$ the space of all pairs $(V,D)$ for $V\in\curly{V}$ and
a K3 surface $D\in|-K_V|$. Then the image of the forgetful map
	\begin{equation}\label{proj}
\pi:(V,D)\in\MM\to D\in\K3(\V)
	\end{equation}
is Zariski open (in particular, dense) in~$\K3(\V)$.
	\end{thrm}
	\begin{remark}
As will be discussed in the next sections, for a number of deformation
classes~$\V$, specific properties of 3-folds in $\V$ ensure that the
map $\pi$ is actually surjective. We do not know if $\pi$ is
surjective for an arbitrary deformation class $\V$ of Fano 3-folds.
	\end{remark}
Theorem~\ref{onto} is proved in the next section.
    \begin{pf}[Proof of Theorem~\ref{ult} assuming Theorem~\ref{onto}]
With the above preparations, the proof will be accomplished almost
entirely via constructions on the K3 lattice~$L$. (The reader
unaccustomed to the K3 lattice may prefer to assume in the first
reading that the Fano 3-folds in $\V_1,\V_2$ have Betti number
$b^2=1$. This would simplify certain technical issues of lattice
arithmetics below and make the principal argument more
transparent while still applicable to some examples.)

Let $D$ be a K3 surface in $\K3\,(\V_1)$ and fix an
isometry between $H^2(D,\ZE)$ and~$L$. This makes $S(\V_1)$ into a
primitive sublattice of~$L$ of rank $b^2(V_1)$, $V_1\in\V_1$. Denote
by $[\kappa_1]$ the image of the induced K\"ahler class of
$D\subset V_1$, so $[\kappa_1]\in S(\V_1)\subset L$ is a primitive
vector of positive square $2n_1-2$ say ($n_1\ge 2$).
By the global Torelli theorem, the complex structure on $D$ is
determined by the image of the complex line $H^{2,0}(D)\subset
H^2(D,\CX)$ in $L\otimes\CX$ under the isometry. The space
$H^{2,0}(D)$ is equivalent to $P(D)=H^2(D,\RE)\cap
(H^{2,0}\oplus H^{0,2})(D)$, an oriented positive-definite real
2-plane in $H^2(D,\RE)$.  For a K3 surface $D$ in $\K3\,(\V_1)$,
the plane $P(D)$ is necessarily orthogonal to the sublattice
$S(\V_1)$ of the Picard lattice of~$D$.
Conversely, a standard result for the K3 surfaces, known as
surjectivity of the period map (\cite{BPV} \S VIII.14), implies that
{\em every} positive 2-plane $P$ in $L\otimes\RE$ orthogonal to
$S(\V_1)$ occurs as $P(D)$ for some $D\in\K3(\V_1)$.

As was explained above, the lattice $S(\V_1)^\bot$ contains $H$ and
therefore for any given $n_2\ge 2$, we can find a primitive vector
$[\kappa_2]\in S(\V_1)^\bot$ of positive length
$([\kappa_2],[\kappa_2])=2n_2-2$. In particular, we may take (we
shall) $2n_2-2$ to be the square of the K\"ahler class of K3 surfaces
in $\K3(\V_2)$. It is not difficult to check that
$[\kappa_1]$,$[\kappa_2]$ can be extended to a basis of~$L$, hence
these two vectors generate a primitive sublattice of~$L$ (this will
become important soon).

If the map $\pi$ of Theorem~\ref{onto} is not surjective for $\V=\V_1$
then one more step in the construction is needed. We shall say that a
K3 surface $D\in\K3(\V_1)$ is {\em $\V_1$-generic} (or just generic
when a confusion is unlikely) if $D$ is in the image of the projection
$\pi$, $D\in\pi(\Mm)$. Now the family of positive 2-planes $P$ through
$[\kappa_2]$ and orthogonal to $S(\V_1)$ defines a real-analytic
subvariety $\K3(\V_1)^\RE$ of {\em real} dimension $20-b^2(V_1)$ in
the moduli space $\K3(\V_1)$ of complex dimension $20-b^2(V_1)$. The
real subvariety $\K3(\V_1)^\RE$ is locally modelled on the projective
space $P^+_1$ of real lines in the positive cone of $L_1\otimes\RE$,
where
$L_1=(S(\V_1)\oplus\langle[\kappa_2]\rangle)^\bot$ is a Lorentzian
sublattice defined by taking the orthogonal complement in~$L$.
On the other hand, local Torelli theorem implies that the complex
moduli space $\K3(\V_1)$ is locally biholomorphic to the period domain
$\{\omega\CX\in \mathbb{P}(S(\V_1)^\bot\otimes\CX):
(\omega,\omega)=0, (\omega,\bar{\omega})>0\}$.
We deduce by examination of the tangent spaces that
$\K3(\V_1)^\RE$ is {\em not} contained in any complex subvariety of
positive codimension in~$\K3(\V_1)$. Thus $\K3(\V_1)^\RE$ cannot lie
in any Zariski closed subset of $\K3(\V_1)$ and has to contain points
defined by generic K3 surfaces. The complement of these generic period
points in $\K3(\V_1)^\RE$ is the common zero set of a finite system of
real-analytic functions. By the identity theorem, this analytic subset
cannot contain open neighbourhoods, so the points in $\K3(\V_1)^\RE$
defined by generic K3 surfaces form a dense open subset.

Now everything in the above constructions concerning $\K3(\V_1)$ can
be repeated for $\K3(\V_2)$, with the obvious change of notation. We
then obtain a sublattice $S(\V_2)$ of $L$ containing a primitive
vector $[\kappa'_2]$ of length $2n_2-2$, corresponding to the K\"ahler
class of a generic projective surface in $\K3(\V_2)$. We also obtain a
primitive vector $[\kappa'_1]\in L$ of positive length $2n_1-2$
orthogonal to $S(\V_2)$ and such that $[\kappa'_1]$ and $[\kappa'_2]$
generate a primitive rank 2 positive-definite sublattice of~$L$
isometric to the previously obtained sublattice generated by
$[\kappa_1]$ and $[\kappa_2]$. By~\cite{BPV}, Theorem~I.2.9(ii), any
two primitive embeddings of a rank~2 even lattice into~$L$ can be
intertwined by an isometry of~$L$, so we may assume without loss of
generality that $[\kappa'_1]=[\kappa_1]$ and $[\kappa'_2]=[\kappa_2]$.

For the final step in the proof, consider first the special case when
the Fano 3-folds in $\V_j$, $j=1,2$, have $b^2=1$, so the sublattice
$S(\V_j)$ is generated by $[\kappa_j]$. Then, by the construction and
by the surjectivity of the period map, any vector, say $[\kappa_K]\in
L\otimes\RE$, of positive length and orthogonal to both $[\kappa_1]$
and $[\kappa_2]$ simultaneously gives us {\em two} K3 surfaces, say
$D_1\in\K3(\V_1)^\RE$ and $D_2\in\K3(\V_2)^\RE$. The positive 2-plane
$P(D_1)$ is spanned by $[\kappa_2]$ and $[\kappa_K]$, whereas $P(D_2)$
is spanned by $[\kappa_1]$ and $[\kappa_K]$. Furthermore, a vector
$[\kappa_K]$ can be chosen so that the two K3 surfaces are generic, as
this just amounts to choosing an element $[\kappa_K]\RE$ in the
intersection of the two dense open subsets of $P^+_1,P^+_2$ identified
above, as in this case $P^+_1=P^+_2 \subseteq
\mathbb{P}(([\kappa_1]\RE\oplus[\kappa_2]\RE)^\bot)\subset 
\mathbb{P}(L\otimes\RE)$.

It remains to multiply the classes $[\kappa_i]$ by constants to define
$[\kappa_I]$,$[\kappa_J]$ such that
$[\kappa_I]^2=[\kappa_J]^2=[\kappa_K]^2$.  Note that this rescaling
changes the K\"ahler classes but does not change the complex
structures of $D_1$ and $D_2$. Thus we have constructed $D_1$ and
$D_2$ satisfying the matching condition, as defined in~\S\ref{K3}.
These K3 surfaces are generic and occur as anticanonical divisors in
Fano 3-folds $V_1$,$V_2$ in $\curly{V}_1$,$\curly{V}_2$ respectively,
by Theorem~\ref{onto}. Therefore, and in view of the remarks earlier
in this section, the construction of a $G_2$-manifold from $V_1$ and
$V_2$ can go ahead. This completes the proof for the Fano 3-folds
having $b^2=1$.

In general, when Fano 3-folds may have $b^2\ge 1$, the choice of
a suitable $[\kappa_K]$ follows the same idea but is slightly
technical. We shall need an auxiliary lemma. 
Recall from the above the definitions of $L_1$ and $P^+_1$.
	\begin{lemma}\label{eichler}
Let $U_1$ be a dense open set in $P^+_1$ and $C^+$ the image of the
positive cone in $\mathbb{P}(S\otimes\RE)$, where
$S\subset L_1$ is a non-degenerate primitive sublattice of signature
$(1,1)$. Then there exists an isometry of~$L$ fixing $S(\V_1) \cup
[\kappa_2]$ and mapping an element of $U_1$ into~$C^+$.
	\end{lemma}
	\begin{pf}
We assume that $C^+\cap U_1$ is empty since otherwise there is
nothing to prove. As the lattice $S$ is indefinite it contains an
isotropic vector $x$ which may be chosen primitive. Let $x'$ be such
that $\{x,x'\}$ is a basis of~$S$. An open set $U_1$ contains a
rationally defined element $e\RE$, with $e\in L_1$ of positive length.
Multiplying $e$ by an integer if necessary, we may write 
$e=ax+bx'+(x,x')b f$, $a,b\in\ZE$, for some $f\in L_1$, $f\neq 0$,
orthogonal to~$S$. As $U_1$ is dense open, we may assume  $b\neq 0$ in
the latter formula. Then Eichler's `elementary transformation' with
respect to $f$ and $x$,
$$
\epsilon_{f,x}: v\in L\to v+(v,f)x-\pol(f,f)(v,x)x-(v,x)f\in L,
$$
provides the required isometry of $L$, as it maps $e$ to the span of
$x$ and~$x'$. As $x$ and $f$ are both orthogonal to $S(\V_1) \cup
[\kappa_2]$, the latter set is fixed by $\epsilon_{f,x}$.
	\end{pf}

Recall that we have two different views on the K3 lattice $L$ (and on
the vectors $[\kappa_1]$,$[\kappa_2]$ in~$L$) depending on whether we
are considering K3 surfaces in $\K3(\V_1)^\RE$ or in $\K3(\V_2)^\RE$.
The sublattices $S(\V_1)$,$S(\V_2)$ together generate a non-degenerate
sublattice of $L$ of signature $(2,t)$, $t\le 18$, so the orthogonal
complement in $L$ of this sublattice has rank at least~2 and contains
a non-degenerate sublattice $S$ of signature (1,1).
Take up one point of view on~$L$. Applying, if necessary,
Lemma~\ref{eichler} to  $S$, with $U_1$ the subset of generic
points in $\K3(\V_1)^\RE$, we may assume that there is a vector
$[\kappa_K]\in S\otimes\RE$ defining a generic K3 surface $D_1$ in
$\K3(\V_1)^\RE\subset\K3(\V_1)$. Moreover the set of the `bad'
elements $[\kappa_K]\RE\in C^+$ giving 
non-generic K3 surfaces in $\K3(\V_1)^\RE$ is an analytic subset,
hence discrete in the 1-dimensional~$C^+$. Note that we are
looking at the structure $(S(\V_1)\cup [\kappa_2])\subset L$, but
the sublattice $S(\V_2)$ is {\em disregarded} for the moment (except
for the vector~$[\kappa_2]$), i.e.\ the action of the isometry given in
Lemma~\ref{eichler} is not applied to $S(\V_2)$. Now shifting to the
other view of~$L$ we consider $S(\V_2)$ (and temporarily disregarding
$S(\V_1)$ except for the $[\kappa_1]$) we apply, if necessary,
Lemma~\ref{eichler} again to the same $S$ but now with $U_2$, $P^+_2$,
$L_2$, respectively defined using the sublattice $S(\V_2)$. This
similarly gives us another dense open subset in~$C^+$, consisting of
those $[\kappa_K]\RE$ 
in the positive cone of $S\otimes\RE$ which define $\V_2$-generic K3
surfaces. Intersecting the two dense open subsets of~$C^+$ we get a
choice of $[\kappa_K]$ that we want. This completes the proof of
Theorem~\ref{ult}.
    \end{pf}

\section{Infinitesimal deformations of the anticanonical K3 divisors}
\label{pf.onto}

In this section we prove Theorem~\ref{onto}: that every K3 surface $D$
in a Zariski open subset in the moduli space $\K3(\V)$ of
lattice-polarized K3 surfaces occurs as an anticanonical divisor in
some Fano 3-fold $V\in\V$. In the case, when the Fano 3-folds in $\V$
have $b^2=1$ (and so $\K3(n)$ is the moduli space of all the
projective K3 surfaces of appropriate degree) this is a known result,
proved in~\cite{CLM} by a different method. The results of that work,
in particular, give a list of (deformation classes of) those Fano
3-folds with $b^2=1$ for which {\em every} K3 surface (of appropriate
degree) occurs as an anticanonical divisor on some member of $\V$,
i.e.\ in our notation, the projection $\pi$ is then surjective.

The argument given below does not use any assumption on Betti numbers
or consideration of the list of Fano 3-folds. The basic idea is that the
space of pairs $\MM$
can be parametrized by a Zariski open set in a (complex) projective
variety, similarly to (an irreducible component of) Chow parameter space
for all the smooth varieties of a given degree in a given projective
space. Cf.\ \S 5.4 in~\cite{CLM1}. Then $\pi$ may be thought of as a
regular map between quasiprojective varieties. By a theorem of
Chevalley, \cite{H}~p.94, the image of $\pi$ is then a `constructible
subset' of $\K3(n)$, which means that it is a finite disjoint union of
sets, each one being the intersection of a Zariski open and a Zariski
closed subset of $\K3(n)$. We have checked that $\K3\,(\V)$ is
connected, so it now suffices to show that the image of $\pi$ has the
right dimension to obtain that $\pi(\MM)$ is Zariski open in
$\K3(\V)$, as claimed.
	\begin{prop}
The dimension of $\pi(\MM)$ is equal to the dimension $20-b^2(V)$
of~$\K3\,(V)$.
	\end{prop}
	\begin{pf}
We shall find the dimension by examining the action of~$\pi$ on the
local first order infinitesimal deformations of a pair $(V,D)$
in~$\MM$. The argument relies on the familiar
Kodaira--Spencer--Kuranishi theory of deformations of the holomorphic
structures~\cite{kodaira}, see also \S2 in \cite{DF} for a related
discussion and further references. We recall, in outline, that the
isomorphism classes of the deformations of a compact complex manifold
$D$, respectively~$V$, define elements in \v{C}ech cohomology group
$H^1T_D$, $H^1T_V$, where $T_D$,$T_V$ denote the sheaves
of holomorphic local vector fields. On the other hand, any
infinitesimal deformation in $H^1T_V$,$H^1T_D$ extends to a genuine
deformation of the complex manifold if the obstruction space $H^2T_V$,
$H^2T_D$, vanishes.  This is the case in the present
situation as, using Serre duality and the Kodaira vanishing theorem
(\cite{GH}, Ch.1), we have $H^2T_V=H^1(\Omega_V^1(-D))=0$ for a
negative line bundle $[-D]$ on a Fano 3-fold and
$H^2T_D=H^0\Omega^1_D=H^{0,1}(D)=0$ on a K3 surface. Notice also that
similarly $h^1T_D=h^{1,1}(D)=20$ equals the dimension of the
moduli space~$\K3$ of all complex K3 surfaces and $H^0T_D=H^{2,1}(D)=0$
verifies that the group of automorphisms of a K3 surface is discrete.

A suitable way to deal with deformations of divisors coming from
deformations of the ambient manifold is discussed in \S1 of~\cite{welters}.
We briefly review the technique adapting it to our situation when the
line bundle is determined by the complex manifold~$V$.
To deal with the pairs $(V,D)$, consider the morphism (more formally,
the complex) of sheaves
$$
d_s: T_V\oplus T_D\to T_V|_D,
$$
where $d_s$ is the sum of the morphisms appearing in
the tangential exact sequence
	\begin{equation}\label{tang}
0\to T_D\to T_V|_D\to N_D\to 0
	\end{equation}
and in the exact sequence
	\begin{equation}\label{k3div}
0\to T_V(-D)\to T_V\to T_V|_D\to 0,
	\end{equation}
obtained by tensoring the structure sequence of $D$ with
$T_V$. Associated to $d_s$ is the double complex of
\v{C}ech cochains with respect to an (appropriately refined) open
covering $\U$ of~$V$,
$$
\begin{CD}
C^0(\U,T_V\oplus\iota_*T_D) @>>> C^1(\U,T_V\oplus\iota_*T_D)
@>>>  \ldots\\
@VV{d_s}V 	@VV{-d_s}V\\
C^0(\U,T_V|_D)    @>>>	C^1(\U,T_V|_D)
@>>>  \ldots .
\end{CD}
$$
The homology of the corresponding single complex is called the
hypercohomology $\HC^*$ of~$d_s$. The hypercohomology
has the same meaning for the deformations of pairs $(V,D)$ as the
cohomology of tangent sheaves does for the deformations of complex
manifolds. In particular, the isomorphism classes of linear
infinitesimal deformations of the pair $(V,D)$ are canonically
parametrized by the first hypercohomology group $\HC^1$ of $d_s$
(cf.~Proposition~1.2 in~\cite{welters}).

The spectral sequence of hypercohomology (see \cite{GH} \S3.5)
of~$d_s$ yields an exact sequence
$$
0\to\HC^1\xra{\beta} H^1T_V\oplus H^1T_D \xra{d_s} H^1(T_V|_D)
\to\HC^2\to 0 ,
$$
where we took account of the vanishing of $H^2T_V$ and $H^2T_D$.
The monomorphism $\beta$ is the `infinitesimal embedding' of the
deformations of the pair $(V,D)$ into the pairs of deformations of $V$
and $D$ alone. We are actually interested in the
composition $\beta_2:\HC^1\to H^1T_D$ of $\beta$ with the forgetful
projection to~$H^1T_D$.

To understand $\beta_2$, first take the cohomology of~\eqref{tang}
to obtain an exact sequence
	\begin{equation}\label{seq}
H^0N_D \to H^1T_D\to H^1(T_V|_D)\to 0.
	\end{equation}
Here we used the isomorphism $N_D=\OH_V(D)|_D$ (the adjunction
formula), then Serre duality on~$D$ and the Kodaira vanishing to claim
$H^1N_D=H^1(\OH_V(D)|_D)= H^1(\OH_V(-D)|_D)=0$. Then
\eqref{seq} tells us that the map  $H^1T_D\to H^1(T_V|_D)$ is
surjective with the kernel given by those infinitesimal deformations
of $D$ which are obtainable by moving $D$ in~$|-K_V|$. It follows that 
$d_s$ is surjective and the obstruction space $\HC^2$ vanishes, hence
{\em every} class in~$\HC^1$ arises from a genuine deformation of the
pair $(V,D)$ and so the image of $\beta_2$ coincides with the image of
the infinitesimal action of $\pi$ on $(V,D)$.

Finally, taking the cohomology of~\eqref{k3div} we can write an exact
sequence
$$
H^1T_V\to H^1(T_V|_D)\xra{\gamma} H^2(T_V(-D))\to 0
$$
(remember $H^2T_V=0$). Using Serre duality we write $H^2(T_V(-D))=
H^1\Omega^1_V= H^{1,1}(V)=H^2(V,\CX)$.
As $\gamma$ is surjective, the image of $H^1T_V$
has (complex) codimension $b_2=b_2(V)$ in $H^1(T_V|_D)$. Then, putting
all the above conclusions together, we calculate that the dimension of
the image of $\beta_2$ is precisely $20-b_2(V)$.
	\end{pf}

\section{New topological types of compact $G_2$-manifolds}
\label{examples}

We next look at the topology of the compact $G_2$-manifolds
constructed from pairs of Fano 3-folds. The notation is as set up
in~\S\ref{fano}.

Corollary~\ref{bFano} and Theorem~\ref{approx}(i) imply that any
$G_2$-manifold constructed from two Fano 3-folds is simply-connected.
We then consider the integral homology groups. As $H_1$ is
determined by the fundamental group and in
view of Poincar\'e duality, we restrict attention to $H_2$
and~$H_3$. Much of the calculation is the linear algebra of
Mayer--Vietoris exact sequences and is fairly straightforward.

Consider first the blow-up $\tV$ of a Fano 3-fold $V$ along a curve
$C=D\cdot D$. The genus of~$C$ is usually referred to
as the genus $g(V)$ of Fano 3-fold and is calculated by the
formula $g(V)=-K_V^3/2+1=\langle c_1(V)^3,[V]\rangle/2 +1$
(e.g.~\cite{Isk77}, Proposition~1.6). Note that $-K_V^3$ is
necessarily positive on a Fano 3-fold~$V$, in particular $g(V)\ge 2$.
Now the exceptional divisor on $\tV$ is a projective bundle of the
Riemann spheres over $C$ and applying standard results on the
cohomology of a blow-up (\cite{GH},~p.605) we find that 
	\begin{equation}
	\label{bup}
H_2(\tV)\cong H_2(V)\oplus\ZE \qquad\text{and}\qquad
H_3(\tV)\cong H_3(V)\oplus\ZE^{2g(V)}.
	\end{equation}

We have already calculated the effect of removing a K3 divisor $D$ from
$\tV=\oW$ on the fundamental group, in the course of proving
Theorem~\ref{line}.
Further application of Mayer--Vietoris theorem to
$\oW=U\cup W$, with a neighbourhood $U$ contracting to~$D$, yields
	\begin{equation}\label{minusD}
H_2(\oW)=H_2(W)\oplus\ZE \qquad\text{and}\qquad
H_3(W)=H_3(\oW)\oplus\ZE^{22-b_2(V)}
	\end{equation}
(as $b_2(V)=b_4(\oW)-1$). For the second isomorphism
in~\eqref{minusD}, recall first that $H_i(Y\times S^1)=H_i(Y)\oplus
H_{i-1}(Y)$ for any manifold~$Y$. Then note that the boundary
homomorphism $H_4(\oW)\to H_3(D\times S^1)=H_2(D)$ evaluates the
intersection cycle with~$D$. The kernel is generated by the cycle
of~$D$ in $H_4(\oW)$, because $D\cdot D=0$ in~$\oW$, but if a divisor
$D'$ is not in the anticanonical class then $D'\cdot D\neq 0$ (this
can be calculated on $V$, away from the blow-up locus~$C$).

Now for the compact real 7-manifold~$M$. Up to a homotopy equivalence,
$$
M\sim(W_1\times S^1)\cup(W_2\times S^1),
\quad\text{such that}\quad
(W_1\times S^1)\cap(W_2\times S^1)\sim D\times S^1\times S^1.
$$
From Mayer--Vietoris theorem and the
knowledge of the homology of K3 surface we can obtain an exact
sequence
$$
H_2(D)\oplus\ZE\xra{\gamma_1}
H_2(W_1)\oplus H_2(W_2)\to H_2(M)\to 0,
$$
as $H_1(W_i)$, $i=1,2$, is trivial. Hence
	\begin{equation}\label{h2}
H_2(M)=\frac{H_2(W_1)\oplus H_2(W_2)}{\gamma_1(H_2(D))}.
	\end{equation}
To understand the homomorphism $\gamma_1$ recall that each of the two
push-forward homomorphisms $H_2(D)\to H_2(V_j)$ is surjective by
Lefschetz hyperplane theorem. From the way we identified the matching
K3 surfaces in the construction of~$M$ we can further deduce that the
image of~$\gamma_1$ contains non-zero elements
$(0,\gamma_1(\mathrm{P.D.}[\kappa_1]))$ and
$(\gamma_1(\mathrm{P.D.}[\kappa_2]),0)$.
The group $H_2(M)$ will have no torsion if both $H_2(W_i)$ are
torsion-free. (For the 
latter, e.g.\ $H^3(V_i)=0$ will suffice.) The dimension of
$\gamma_1(H_2(D))$, and hence the Betti number $b_2(M)$, is not in
general uniquely determined by $V_1,V_2$ because it depends on the
rank of the sublattice of $H^2(D)$ generated by $S(\V_1)$ and
$S(\V_2)$. In any event, we have
	\begin{equation}\label{range}
b_2(M)\le\min\{b_2(V_1),b_2(V_2)\}-1.
	\end{equation}
On the other hand, we shall see an example below when all the 
values of $b_2(M)$ allowed by~\eqref{range} actually occur.
As any Fano \mbox{3-fold} $V_j$ has $1\le b_2(V_j)\le 10$, any
$G_2$-manifold $M$ constructed from a pair of Fano 3-folds necessarily
satisfies
$$
0\le b_2(M)\le 9.
$$

The group $H_3(M)$ is found by similar, although slightly lengthy
computations, taking account of the previous steps. We omit the
details. Note, in particular, that the inclusion homomorphism of
$H_3(D\times S^1\times S^1)$ is equivalent to the direct sum of two maps
$H_2(D)\to H_2(W_{3-i})\oplus H_3(W_i)$, $i=1,2$, and these maps miss
the subspaces $H_3(\oW_i)\subset H_3(W_i)$. The dimension counting once
again depends on the rank of the sublattice generated by $S(\V_1)$ and
$S(\V_2)$ and yields (substituting $\tV_j=\oW_j$)
	\begin{equation}\label{b3M}
b_3(M)+b_2(M)=b_3(\tV_1)+b_3(\tV_2)+23,
	\end{equation}
for any $G_2$-manifold $M$ constructed from $V_1$ and $V_2$.

In the case when one of the two Fano 3-folds has $b_2=1$ we can say more.
	\begin{thrm}\label{b3}
Suppose that $M$ is a compact 7-manifold with holonomy $G_2$
constructed from Fano 3-folds $V_1$ and $V_2$, where $b^2(V_1)=1$.
Then
$$
b_2(M)=0\qquad\text{and}\qquad
b_3(M)=b_3(V_1)-K^3_{V_1}+b_3(V_2)-K^3_{V_2}+27.
$$
Further, the diffeomorphism type of $M$ (as a smooth real 7-manifold)
is independent of the choice of $V_1$ and $V_2$ in their deformation
classes.
	\end{thrm}
	\begin{pf}
The assertions, except for the last sentence, are immediate from the
above topological considerations.  The diffeomorphism type of a smooth
real 6-manifold underlying $W$ is not affected by moving an
anticanonical K3 divisor or deforming the complex structure of a
Fano 3-fold. Thus the only possible ambiguity in defining a smooth
7-manifold~$M$ is in the choice of the matching diffeomorphism of the
K3 divisors (recall \S\ref{7mfd}). Let a piece of~$M$, diffeomorphic
to $W_1\times S^1$, be fixed and consider two choices of matching K3
surfaces $D_1\to D_2$ and $D_1\to D'_2$. Respectively, there are two
choices $W_2\times S^1$ or $W'_2\times S^1$ for attaching the other
piece; denote the resulting 7-manifolds by $M$ and~$M'$. It then
remains to demonstrate that the identity map of $W_1\times S^1$ can be
extended to a diffeomorphism $M\to M'$. The
$\mathrm{id}_{W_1\times S^1}$ readily extends to the diffeomorphism
between the cylindrical ends $D\times S^1\times S^1\times\RE$ of
$W_2\times S^1$ and $W'_2\times S^1$ which is determined by a
diffeomorphism of the real 4-manifold~$D$ underlying K3 surfaces
$D_2,D'_2$. We claim that this latter diffeomorphism is isotopic (can
be connected by a smooth path  of automorphisms) to the identity map
of~$D$. Assuming this claim, it is easy to see that $M$ is diffeomorphic
to~$M'$. To verify the claim, we refer back to the lattice
considerations and the notation in the proof of Theorem~\ref{ult}. A
choice of $D_2$ or $D'_2$ corresponds, in the first place, to a choice
of positive lattice vector $[\kappa_2]$, orthogonal to a given
$[\kappa_1]$, but all the choices of $[\kappa_2]$ are equivalent up to
an isometry of the K3 lattice~$L$. (It is at this point that the
condition $b^2(V_1)=1$ is used.) Then the ambiguity comes down to a
choice of positive $[\kappa_K]\in L\otimes\RE$. But any two such
choices are connected by a path of positive vectors. Then Torelli
theorem and surjectivity of the period map translate this into a path
of complex structures from $D_2$ to $D'_2$ and induce the required
path of the diffeomorphisms of the underlying real manifold~$D$.
	\end{pf}

The $G_2$-manifolds with vanishing or finite second homology group
give a new series. Previously, only one topological type of
\mbox{$G_2$-manifold} with $b_2=0$ was known (it has \mbox{$b_3=215$),}
obtained by resolving singularities of the quotient of 7-torus by
$\ZE_2^2$, see~\cite{joyce-book}, \S 12.5.5.

Here are examples of new $G_2$-manifolds constructed from some
well-known Fano manifolds and treated rather more explicitly than in
the general results in the previous sections.
	\begin{example}[Rigid Fano 3-folds]\label{rigid}
$V=\CP^3$ is a Fano 3-fold and any smooth quartic is a K3 surface
in the anticanonical class, so in this example $\pi$ maps
$\curly{P}(CP^3)$ {\em onto} $\K3(3)$. More generally, we remark that
in the case when the deformation class $\V$ consists of $V$ alone the
image of the projection $\pi:\MM\to\K3(\V)$ is a {\em closed} subset
in $\K3(\V)$ and hence $\pi$ is surjective. The closedness is not
difficult to see, for in this case $\MM$ is identified with a Zariski
open subset in the projective space parametrizing the anticanonical
linear system $|-K_V|$ and $\pi$ becomes a tautological map to the
isomorphism class of a divisor.

By Theorem~\ref{ult}, a compact $G_2$-manifold, say $M_0$, can be
constructed from two copies of $\CX P^3$, taking a pair of matching
quartics $D_1,D_2$. As the cycle of any $D=D_i$ is 4 times 
the hyperplane cycle, then $-K_V^3=\pol D\cdot D\cdot D=64$.
The resulting $G_2$-manifold has
$H_2(M_0)=0$ by~\eqref{h2}, and $b_3(M_0)=155$ by Theorem~\ref{b3}. 

A Fano 3-fold $V'=\CP^2\times\CP^1$ has $H^2(V')=\ZE\oplus\ZE$, so
in this case we encounter anticanonical divisors with a rank~2
polarization lattice $S(V')$. The anticanonical divisors on $V'$ are
defined by the polyhomogeneous polynomials of bidegree $(3,2)$ in the
respective homogeneous variables. Calculation in the cohomology ring
gives $-K^3_{V'}=54$. Taking the intersections of 2-cycles defined on
$D'\in|-K_{V'}|$ by $\CP^2\times\mathrm{pt}$ and $\ell\times\CP^1$,
where $\ell$ denotes a projective line in $\CP^2$, we find
$$
S(V')=\begin{pmatrix} 0&3\\ 3&2 \end{pmatrix}.
$$
The induced K\"ahler class on~$D$ is Poincar\'e dual to the cycle
expressed by the vector $(1;1)$, with respect to the basis of $S(V')$.
So $D'$ may be realized as an $S(V')$-polarized octic K3 surface,
indeed it embeds in $\CP^5$ by Segre embedding of~$V'$.
It follows that a $G_2$-manifold
$M'$ constructed from $\CP^3$ and $\CP^2\times\CP^1$ has
$H_2(M')=0$ (the vanishing of $H^3(\CP^3)$ and
$H^3(\CP^2\times\CP^1)$ implies that there is no torsion),
and $b_3(M')=145$, by Theorem~\ref{b3}.

For the construction of a $G_2$-manifold from two copies of
$\CP^2\times\CP^1$ one needs to take account of the rank of the
span of two copies of the polarization lattice $S(V')$ embedded in the
K3 lattice (recall the proof of Theorem~\ref{ult}). There are two
possibilities to consider, as the rank can be 3 or 4.

It is convenient to use a basis $\{f_1=(1;1),f_2=(1;0)\}$ of $S(V')$
to keep track of the K\"ahler class~$f_1$ in the embeddings.
The intersection numbers in the new basis are $(f_1,f_1)=8$,
$(f_1,f_2)=3$, $(f_2,f_2)=0$.
Let $\{e_j,e'_j\}$, $j=1,2,3$, be the standard basis of the $j$-th
hyperbolic plane summand~$H$ in the K3 lattice, so
$(e_j,e_j)=(e'_j,e'_j)=0$ and $(e_j,e'_j)=1$. The primitive embedding
of two copies of $S(V')$ defined by $f_1\to [\kappa_1]=e_1+4e'_1$, and
$f_1\to [\kappa_2]=e_2+4e'_2$, with $f_2\to
e_1-e'_1+e_2-e'_2+e_3+2e'_3$ for both copies, have the property that the
span of the two images of $S(V')$ has rank~3. Denote by $M_1$ the
corresponding $G_2$-manifold. Then $\pi_1(M_1)=1$, $H_2(M_1)=\ZE$, and
$b_3(M_1)=134$. The previously published examples of $G_2$-manifolds
with $b_2=1$, constructed by a different method, have either $b_3=142$
or $b_3=186$ (\cite{joyce-book} \S\S 12.5.5 and 12.7.4), so the
$G_2$-manifold $M_1$ is new.

It is also possible to embed two copies of $S(V')$ in the K3 lattice so
that together they generate a lattice of rank~4. This constructs a
$G_2$-manifold $M_2$ which is again simply-connected, but now $H_2(M_2)$
vanishes and $b_3(M_2)=135$. Details are left to the reader.
	\end{example}
	\begin{example}[Branched double covers]
One such Fano 3-fold is the double cover $X_2$ of $\CP^3$ branched
over a smooth quartic surface~$D$. For a double cover over a divisor,
$b_2(X_2)=b_2(\CP^3)$ (note that by Lefschetz--Bott any 2-cycle on a
3-fold may be assumed to be in the branch locus). The Euler
characteristic is $\chi(X_2)=2\chi(\CP^3)-\chi(D)=-16$, hence
$b_3(X_2)=20$. The divisor $D$ lifts isomorphically to a divisor $D'$
on~$X_2$ and calculation in cohomology (cf. \cite{GH}, \S4.4.) shows
$-K_{X_2}=p^*[-D]-D'=2[D']-[D']=[D']$, thus $D'$ is an
anticanonical divisor on $X_2$ (and so every smooth quartic occurs in
the anticanonical class of some~$X_2$). Further, $D'$ has degree two
in $X_2$, so $-K_{X_2}^3=8\cdot 2$, the extra factor 2 coming from the
degree of the covering. Then the topology of a 
$G_2$-manifold $M$ constructed from $X_2$ and any Fano of the previous
example (or $X_2$ again) is obtained by Theorem~\ref{b3}, e.g.\ if
$M$ is constructed from $X_2$ and $\CP^2\times \CP^1$ then
$b_2(M)=0$, $b_3(M)=117$.
	\end{example}
	\begin{example}[Fano 3-folds as divisors]\label{F.div}
A quadric hypersurface $Q\subset\CX P^4$ is a Fano 3-fold with
$b_2(Q)=1$, $b_3(Q)=0$. The value of $b_2$ is implied by Lefschetz
hyperplane theorem and then $b_3$ is again determined by the Euler
characteristic. For the latter, recall that we may express the full
Chern class of a projective hypersurface (and more generally, of a
complete intersection of hypersurfaces) knowing that the degree of its
normal line bundle equals the degree of hypersurface by the adjunction
formula. Then, writing $x$ for the pull-back of the generator of the
cohomology ring of the ambient $\CP^4$ we find
$c(Q)=(1+x)^5(1+2x)^{-1}=1+3x+4x^2+2x^3$, hence
$\chi(Q)=\langle 2x^3,[Q]\rangle=4$. The computation $c_1(Q)=3x$
also tells us that a K3 surface $D$ in $|-K_{Q}|$ is obtained as
complete intersection of $Q$ with a cubic hypersurface
$X_3\subset\CP^4$, so $\deg(D)=6$. Conversely, for any given sextic K3
surface $D$ in $\CP^4$, consider the linear systems of quadric
and cubic 3-folds containing~$D$. We can find {\em non-singular} 
quadric and cubic 3-folds in these linear systems by Bertini's
theorem~\cite{GH}. Note that this proves Theorem~\ref{onto} for the
deformation class of $Q$, including the surjectivity of
the projection~$\pi$. Now $-K^3_{Q}=-\chi(C)$ and notice that a curve
$C=D\cdot D$ is obtained by intersecting $D$ with a generic cubic
hypersurface, so $\deg(C)=18$ and
$c(C)=(1+x)^5(1+2x)^{-1}(1+3x)^{-2}=1-3x$, giving
$-K^3_{Q}=54$.

A cubic $X_3\subset\CX P^4$ is also a Fano 3-fold and similar
arguments show that every sextic K3 surface $D$ occurs in $|-K_{X_3}|$
for some smooth $X_3$, by intersecting with a suitable quadric
hypersurface. We calculate, by the above method, that $b_2(X_3)=1$,
$b_3(X_3)=10$, and $-K^3_{X_3}=24$.
	\end{example}
The calculations are simplified for Fano 3-folds $V$ whose
anticanonical sheaf is very ample, which means that the anticanonical
linear system of $V$ defines an {\em embedding} making $V$ into a
smooth algebraic variety in $\CP^{g+1}$, where $g$ is the genus of $V$.
Then the anticanonical bundle
of $V$ is the restriction of $\OH(1)$ and any K3 surface $D\in |-K_V|$
is a (generic) hyperplane section of~$V$. Further, in this case,
$$
-K^3_V=2g-2=\deg V=\deg D,
$$
so $D\subset\CP^g$ is an algebraic K3 surface of degree $2g-2$.
Of course, a curve $C=D\cdot D$ is just the intersection of~$V$ with
two generic hyperplanes.

A Fano 3-fold is called prime if the class $c_1(V)$ is primitive and
generates the cohomology lattice $H^2(V)$ (then necessarily 
$b^2(V)=1$).
By Theorem~\ref{b3} the
two prime Fanos determine $M$ uniquely up to a diffeomorphism.
	\begin{example}[Complete intersections]\label{ci}
A smooth complete intersection $V=X_8$ of three quadric hypersurfaces
in~$\CP^6$ is a prime Fano 3-fold of degree~8 and genus~5. The prime
property can be seen by first calculating the full Chern class (cf.\
Example~\ref{F.div}) $c(X_8)=(1+x)^7(1+2x)^{-3}=1+x+3x^2-3x^3$. It
follows that $\chi(X_8)=-24$ and as $b_2(X_8)=1$ by Lefschetz, then
$b_3(X_8)=28$.  A hyperplane section $D\in|-K_{X_8}|$ is genericly a
smooth octic K3 surface in $\CP^5$. Applying Bertini's theorem, we see
that every octic K3 surface occurs in the anticanonical divisor on
some (smooth)~$X_8$.

Another similar example of prime Fano 3-fold $V=X_6$ is given by a
smooth complete intersection of a quadric and a cubic in $\CX P^5$. It
has genus~4 and degree~6, with $b_2(V)=1$, $b_3(V)=40$.
	\end{example}
	\begin{example}[Prime Fano 3-fold of genus~12]
This prime Fano 3-fold $X_{22}$, a smooth variety  of degree 22 in
$\CP^{13}$, was originally found by
Iskovskikh~\cite{Isk77} and further investigated 
by Mukai. See \cite{mukai} and references therein for various
descriptions of~$X_{22}$. A 3-fold $X_{22}$ is not isomorphic to a
branched double cover of any homogeneous space or a complete
intersection of hypersurfaces in a homogeneous space. It is known that
$X_{22}$ is a compactification of 
$\CX^3$ and has the homology of $\CP^3$, but admits a 6-dimensional
family of deformations. By the remarks before Example~\ref{ci}, an
anticanonical divisor on any $X_{22}$ defines a point in $\K3(22)$.
In the case of $V=X_{22}$ Theorem~\ref{onto} reproves Theorem~7
of~\cite{CLM1}. As $-K^3_{X_{22}}=22$ Theorem~\ref{b3} shows that
using $X_{22}$ instead of $\CP^3$ one obtains $G_2$-manifolds with
`smaller' homology and fundamental group.
	\end{example}

The $G_2$-manifolds constructed from any two of the Fano 3-folds
described above satisfy $71\le b_3(M)\le 155$, by
Theorem~\ref{b3}. The upper bound is attained by constructing $M$
from a pair of $\CP^3$'s and the lower bound from a pair of prime
$X_{22}$'s.

\end{document}